\documentclass[oneside,a4paper,12pt]{article}
\usepackage{amsmath}
\usepackage{amssymb}
\usepackage{graphicx}
\usepackage{natbib}
\newtheorem{Definition}{\bf \large Definition}[section]
\newtheorem{Theorem}{\bf \large Theorem}[section]
\newtheorem{Lemma}[Theorem]{\bf \large Lemma}
\newtheorem{Corollary}[Theorem]{\bf \large Corollary}
\newtheorem{Remark}{\bf \large Remark}[section]
\renewcommand{\theequation}{\thesection.\arabic{equation}}
\newcommand{\proof}{\noindent \textbf{Proof.} \hspace {0.3cm}}

\newcommand{\pr}{\mathrm{P}}
\newcommand{\E}{\mathrm{E}}
\newcommand{\dif}{\mathrm{d}}

\begin{document}
\title{Optimal financing and dividend distribution in a general diffusion model with regime switching \footnotetext{ Jinxia Zhu, School
      of Risk and Actuarial Studies, The University of New South Wales
      Kensington Campus, NSW 2052, Australia;
      E-mail:jinxia.zhu@unsw.edu.au}
      \footnotetext{Hailiang Yang, Department of Statistics and Actuarial Science, The University of Hong Kong, Pokfulam Road, Hong Kong;
E-mail: hlyang@hku.hk } }
\author{{\sc Jinxia Zhu}\\
 {\it The University of New South Wales,
  Australia}\\{\sc Hailiang Yang}\\
   {\it The University of Hong Kong, Hong Kong}}

\date{}
\maketitle

\begin{abstract}
We study the optimal financing and dividend distribution problem with restricted dividend rates in a diffusion type surplus model where the drift and volatility coefficients are general functions of the level of surplus and the
external environment regime. The  environment regime is modeled by a Markov process. Both capital injections and dividend payments incur expenses.  The objective is to maximize the expectation of the total discounted dividends minus the total cost of capital injections. We prove that it is optimal to inject capitals only when the surplus tends to fall below zero and to pay out dividends at the maximal rate when the surplus is at or above the threshold dependent on the environment regime.
\end{abstract}


\section{Introduction }
The optimal dividend strategy problem has gained extensive attention. In the diffusion setting, many works concerning dividend optimization use the Brownian motion model for  the underlying cashflow process.   \cite{Bauerle2004} extends the basic model by assuming that the drift coefficient is a  linear function of the level of cashflow and  \cite{CadenillasSarkarZapatero2007} uses the mean-reverting model  and solves the optimization problem. \cite{HojgaardTaksar2001} considers the optimization problem under the model where the drift coefficient is proportional to the level of cashflow and the diffusion coefficient is proportional to the square root of the  cashflow level.  \cite{ShreveLehoczkyGaver1984},
  \cite{Paulsen2008}, \cite{Zhu2014b} and some references therein address the optimization problems for the general diffusion model where the drift and diffusion coefficients are general functions of the  cashflow level.

An interesting and different direction of  extension is to include the impact of the changing external  environments/conditions (for example, macroeconomic  conditions and weather conditions) into modeling of the cashflows.  A continuous time Markov chain can be used  to model the state of the external environment condition, of which the use is  supported by observation in financial markets.
The optimal dividend problem with  regular control
for Markov-modulated risk processes has been investigated under a verity of assumptions.  \citet{SotomayorCadenillas2011}  solves the dividend optimization problem for a Markov-modulated Brownian motion model with both the drift and  diffusion coefficients modulated
by a two-state  Markov Chain.   \cite{Zhu2014c} solves the problem for the  Brownian motion model modulated by a multiple state Markov chain.

The optimality results in all the above works imply that distributing dividends according to the optimal strategy leads almost surely to ruin. \cite{DicksonWaters2004} proposes to include capital injections (financing) to prevent the surplus becomes negative  and therefore prevent ruin. Under the Brownian motion, \cite{LokkaZervos2008} investigates the optimal dividend and financing problem, and \cite{HeLiang2008} studies  the problem with risk exposure control through control of reinsurance rate.
The optimality problem with control in both capital injections and dividend distribution in a Cram\'er-Lundberg model is addressed  in \cite{ScheerSchmidli2011}. \cite{YaoYangWang2011} solves the  problem for dual model with  transaction costs.

The purpose of this paper is to investigate optimal financing and dividend distribution problem with restricted dividend rates in a general diffusion model with regime switching. Under the model,  the drift and volatility coefficients are general functions of the level of surplus and the
external environment regime, which is modeled by a Markov process. Similar to the ``reflection problem", the company can control the financing /capital injections  process (a deposit process) and the dividend distribution process (a ``withdrawal" process). Both capital injections and dividend payments will incur transaction costs.
 Sufficient capital injections must be made to keep the controlled surplus
process nonnegative and the dividend payment rate is capped. This paper can be considered as an extension of the existing works on the dividend optimization problem with restricted dividend rates for the diffusion models with or without regime switching. The model considered is  more general as it assumes that 1. the drift and volatility are general functions of the cashflows; and  2. the model risk parameters (including drift, volatility and discount rates) are dependent on the external environment regime.

The rest of the paper is organized as follows. We formulate the optimization problem in Section 2. An auxiliary problem is introduced and solved in Section 3. Section 4 presents the optimality results.  A conclusion is provided in Section 5. Proofs  are relegated to Appendix.

\section{Problem Formulation}
Consider a probability  space
$(\Omega,\mathcal{F},\pr)$.  Let $\{W_t;t\ge 0\}$ and $\{\xi_t;t\ge 0\}$ be respectively a standard Brownian motion and  a Markov chain with the finite state space $\mathcal{S}$ and the transition intensity matrix $Q=(q_{ij})_{i,j\in\mathcal{S}}$. The two stochastic processes $\{W_t;t\ge 0\}$ and $\{\xi_t;t\ge 0\}$ are independent.  We use  $\{\mathcal{F}_t;t\ge 0\}$ to denote the minimal complete $\sigma$-field generated by the stochastic process $\{(W_t,\xi_t);t\ge 0\}$. Let $X_t$ denote the surplus at time $t$ of  a firm in absence of financing and dividend distribution. Assume that $X_0$ is $\mathcal{F}_0$ measurable and that $X_t$ follows the  dynamics, $\dif X_t = \mu(X_{t-},\xi_{t-})\dif t+\sigma(X_{t-},\xi_{s-})\dif W_t$ for $t\ge 0$,
where  the functions  $\mu(\cdot,j)$   and $\sigma(\cdot,j)$ are Lipschitz continuous,  differentiable and grow at most linearly on $[0,\infty)$ with $\mu(0,u)\ge 0$.  Furthermore, the function $\mu(\cdot,j)$ is concave and the function $\sigma(\cdot,j)$ is positive and non-vanishing.

The firm must have nonnegative assets in order to continue its business. If necessary, the firm needs to raise money from the market. For each dollar of money raised, it includes  $c$ dollars of transaction cost and hence leads to an increase of $1-c$ dollars in the surplus  through capital injection. Let $C_t$ denote the cumulative amount of   capital injections up to time $t$. Then the total cost for capital injections  up to time $t$ is  $\frac{C_t}{1-c}$.  The company can distribute part of its assets to the shareholders as dividends. For each dollar of dividends received by the shareholders, there will be $d$ dollars of cost incurred to them. Let $D_t$ denote the cumulative amount of dividends paid out by the company up to time $t$. Then the total amount of dividends received by the shareholders up to time $t$ is $\frac{D_t}{1+d}$. We consider the case where the dividend distribution rate is restricted. Let the random variable $l_s$ denote the  dividend  payment rate at time $s$ with the restriction  $0\le l_s\le \bar{l}$ where $\bar{l}(>0)$ is constant.
Then
$D_t=\int_0^tl_s\dif s.$ Both $C_t$ and $D_t$ are controlled by the company's decision makers. Define $\pi=\{(C_t,D_t);t\ge 0\}$. We call $\pi$ a control strategy.

  Taking financing and dividend distribution into consideration,
the dynamics of the (controlled) surplus process with the strategy $\pi$ becomes
\begin{eqnarray}
  \label{eq:2}
  \dif X_t^\pi
  &=&(\mu(X_{t-}^\pi,\xi_{t-})-l_t)\dif t+
  \sigma(X_{t-}^\pi,\xi_{t-})\dif W_t+ \dif C_t,\ t\ge 0.
\end{eqnarray}
Define
$\pr_{(x,i)}\left(\ \cdot\
    \right)=\pr\left(\ \cdot\ |X_0=x,\xi_0=i
    \right),$ $\E_{(x,i)}\left[\ \cdot\
    \right]=\E\left[\ \cdot\ |X_0=x, \xi_0=i
    \right],$
    $\pr_{i}\left(\ \cdot\
    \right)=\pr\left(\ \cdot\ |\xi_0=i
    \right),$ and $\E_{i}\left[\ \cdot\
    \right]=\E\left[\ \cdot\ | \xi_0=i
    \right].$
The performance of  a control strategy $\pi$ is measured by its return function  defined as follows:
 \begin{align}
  R_\pi(x,i)=\E_{(x,i)}\left[\int_0^{\infty} e^{-\Lambda_t}\frac{l_t}{1+d}\dif t-\int_0^{\infty} e^{-\Lambda_t}\frac1{1-c}\dif C_t
    \right],\ x\ge 0,i\in \mathcal{S},\label{jf1}
\end{align}
where $\Lambda_t=\int_0^t\delta_{\xi_s}\dif s$ with $\delta_{\xi_s}$ representing the  force of discount at time $s$. Assume $\delta_i>0$, $i\in\mathcal{S}$.

A strategy $\pi=\{(C_t,D_t);t\ge 0\}$ is said to be \textit{admissible} if (i) both $\{C_t;t\ge 0\}$ and $\{D_t;t\ge 0\}$ are nonnegative, increasing, c\`adl\`ag, and $\{\mathcal{F}_t;t\ge 0\}$-adapted processes, (ii) there exists an $\{\mathcal{F}_t;t\ge 0\}$-adapted process $\{l_t;t\ge 0\}$ with $l_t\in[0,\bar{l}]$ such that  $D_t=\int_0^tl_s\dif s$ and (iii) $X^\pi_t\ge 0$ for all $t>0$. We use $\Pi$ to denote the class of admissible strategies.

Since $\{C_t;t\ge 0\}$ is right continuous and increasing, we have the following decomposition: $C_t=\tilde{C}_t+C_{t}-C_{t-}$, where $\{\tilde{C}_t;t\ge 0\}$ represents the continuous part of $\{C_t;t\ge 0\}$.

  For convenience, we use $X$,  $X^\pi$, $\xi$ and $(X^\pi,\xi)$ to denote the stochastic processes
 $\{X_t;t\ge 0\}$,  $\{X_t^\pi;t\ge 0\}$, $\{\xi_t;t\ge 0\}$ and $\{(X_t^\pi,\xi_t);t\ge 0\}$,
respectively.
Note that for any admissible strategy $\pi$, the stochastic process $X^\pi$ is right-continuous and adapted to the filtration $\{\mathcal{F}_t;t\ge 0\}$.

The objective of this paper is to study  the maximal return function (value function):
\begin{align}
  V(x,i)=\sup_{\pi\in \Pi}R_\pi (x,i),\label{13613-1}
\end{align} and to identify the associated optimal admissible strategy, if any.
Following the standard argument in  stochastic control theory
\citep[e.g.][]{FlemingSoner1993}, we can show that the value function
fulfils the following dynamic programming principle:  for any stopping time $\tau$,
\begin{align}
  V(x,i)=\sup_{\pi\in\Pi}\E_{(x,i)}\Big[\int_0^{ \tau} \frac{l_te^{-\Lambda_t}}{1+d}\dif t-\int_0^{ \tau} \frac{ e^{-\Lambda_t}}{1-c}\dif C_t
  + e^{-\Lambda_{ \tau}} V(X^{\pi}_{  \tau},\xi^\pi_{
    \tau})\Big ].\label{60214-3}
    \end{align}

\section{An Auxiliary Optimization Problem}\label{sec3}
Motivated by \cite{JiangPistorius2012}, which  introduces an auxiliary problem where the objective functional is modified in a way such that only the ``returns" over the time period from the beginning up to the first regime switching are included plus a terminal value at the moment of the first regime switching, we start with a similar auxiliary problem first. The optimality results of this  problem will play an essential role in solving the original optimization problem.

Throughout the paper, we define $\underline{\delta}=\min_{j\in \mathcal{S}}\delta_{j}$, $q_i=-q_{ii}$, and $\sigma_1=\inf\{t>0: \xi_t\neq \xi_0\}$. Here, $\sigma_1$ is the first transition time of the Markov process $\xi$. For any function $g:\mathbb{R}^+\times \mathcal{S}\rightarrow \mathbb{R}^+$, we use $g^\prime(\cdot)$ and $g^{\prime\prime}(\cdot)$ to denote the first order and second order derivatives, respectively, with respect to
the first argument. We start with introducing two special classes of functions.
\begin{Definition}
(i) Let
$\mathcal{C}$ denote  the class of functions $g: \mathbb{R}^+\times \mathcal{S}\rightarrow \mathbb{R}$ such that for each $j\in\mathcal{S}$, $g(\cdot,j)$ is  nondecreasing and  $g(\cdot,j)\le \frac{\bar{l}}{\underline{\delta}(1+d)}$.
(ii) Let
$\mathcal{D}$ denote  the class of functions $g\in\mathcal{C}$ such that for each $j\in\mathcal{S}$, $g(\cdot,j)$ is concave and $\frac{g(x,j)-g(y,j)}{x-y}\le \frac1{1-c}$ for $0\le x<y$. (iii) Define the distance $||\cdot||$ by
$
||f-g||=\max_{x\ge 0,i\in \mathcal{S}}|f(x,i)-g(x,i)|\ \mbox{ for $f,g\in\mathcal{D}$}.$
\end{Definition}
\begin{Lemma} \label{complete}The metric space $(\mathcal{D},||\cdot||)$ is complete.
\end{Lemma}

Define a modified return function and the associated optimal return function by
\begin{align}
  R_{f,\pi}(x,i)
  =&\E_{(x,i)}\bigg[\int_0^{ \sigma_1} \frac{l_te^{-\Lambda_t}}{1+d}\dif t-\int_0^{ \sigma_1} \frac{e^{-\Lambda_t}}{1-c}\dif C_t
 +e^{-\Lambda_{\sigma_1}}f(X^\pi_{\sigma_1}, \xi_{\sigma_1})
    \bigg],\ x\ge 0,i\in\mathcal{S},\label{4214-2}\\
    V_f(x,i)=&\sup_{\pi\in \Pi}R_{f,\pi}(x,i),\ x\ge 0,i\in\mathcal{S}.\label{17713-4}
\end{align}
\begin{Lemma}\label{remff9}
For any $f\in\mathcal{C}$,
$V,V_f\in\mathcal{C}$ .
 \end{Lemma}
Notice that the un-controlled  process $(X,\xi)$ is a Markov process.
For any  $f\in\mathcal{C}$ and any $i\in\mathcal{S}$, the following Hamilton-Jacobi-Bellman (HJB) equation for the modified value function $V_f(\cdot,i)$ can be obtained  by using standard arguments in stochastic control: for $x\ge 0$\\
{\small$
\max\big\{\max_{l\in[0,\bar{l}]}\left(\frac{\sigma^2(x,i)}{2}V_f^{\prime\prime}(x,i)+
\mu(x,i)V_f^\prime(x,i)-\delta_i V_f(x,i)+l\left(\frac1{1+d}-V^\prime_f(x,i)\right)\right),V^\prime_f(x,i)-\frac1{1-c}\big\}=0$.}

Now we define a special class of admissible strategies, which has been shown in the literature to contain the optimal strategy for the original optimization problem if there is 1 regime only.
Since the return function of the modified optimization  includes the dividends and capital injections in the first regime only, this problem can be considered as a problem to maximize the returns up to an independent exponential time for a risk model with 1 regime. It is worth studying the special class of strategies mentioned above to see whether the optimal strategy of the modified problem falls into this class as well.

\begin{Definition}
For any $b\ge 0$, define the strategy $\pi^{0,b}=\{(C^{0,b}_t, D^{0,b}_t);t\ge 0\}$ in the way such that the company pays dividends at the maximal rate $\bar{l}$ when  the surplus equals or exceeds $b$, pays no dividends when the surplus is below $b$ and the company injects capital to maintain the surplus  at level $0$ whenever the surplus tends to go below $0$ without capital injections.
\end{Definition}

We now investigate whether a strategy $\pi^{0,b}$ with an appropriate value for $b$ is optimal or  not  for the modified optimization problem. We start with studying the associated return functions.
For convenience, we  write  $X^{0,b}=X^{\pi^{0,b}}$ throughout the rest of the paper.
\begin{Remark}\label{rem1}
 (i) It is not hard to see that $\pi^{0,b}$ is admissible and that both $\pi^{0,b}$ and $X^{0,b}$ are Markov processes. (ii) For any  function $f\in\mathcal{C}$  and  any $i\in\mathcal{S}$, by applying the comparison theorem used to prove the non-decreasing property of $V(\cdot,i)$ and $V_{f}(\cdot,i)$  in Lemma \ref{remff9} we can show that the function $R_{f,\pi^{0,b}}(\cdot,i)$ is non-decreasing on $[0,\infty)$ as well.
\end{Remark}

%
%
%
For any $f\in\mathcal{C}$, $i\in\mathcal{S}$ and $b\ge 0$, define the operator $\mathcal{A}_{f,i,b}$ by
\begin{eqnarray}
\mathcal{A}_{f,i,b}\ g(x)=\frac{\sigma^2(x,i)}{2}g^{\prime\prime}(x)+
(\mu(x,i)-\bar{l})g^\prime(x)-(\delta_i+q_i) g(x)+\frac{\bar{l}}{1+d}+\sum_{j\neq i}q_{ij}f(x,j)=0.\label{agenerator}
\end{eqnarray}

The following conditions will be required for the main theorems.

\noindent \textbf{Condition 1}: The functions $\mu(\cdot,i)$ and $\sigma(\cdot,i)$ are the ones such that   for any  given function $f\in\mathcal{D}$   and any given $i\in\mathcal{S}$, the ordinary differential equation $\mathcal{A}_{f,i,b}\ g(x)=0$
with any finite initial value at $x=0$  has a  bounded solution over $(0,\infty)$.

A sufficient condition for Condition 1 to hold is that both the functions
$\mu(\cdot,i)$ and $\sigma(\cdot,i)$ are  bounded on $[0,\infty)$ (see Theorem 5.4.2 in \cite{Krylov1996}).  However, this is far away from necessary.  For example,  when $\mu(\cdot,i)$  is a  linear function with positive slope and $\sigma(\cdot,i)$ is a constant Condition 1 also holds (see section 4.4 of \cite{Zhu2014b}).

\noindent \textbf{Condition 2}:
 $\mu^\prime(x,i)\le \delta_i$ for all $x\ge 0$ and $i\in\mathcal{S}$.

Define  for any  function $f \in\mathcal{C}$ and $i\in\mathcal{S}$, \begin{eqnarray}A_{f,i}
=\frac{\bar{l}/(1+d)+\sum_{j\neq i}q_{ij}f(\infty,j)}{q_i+\delta_i}.\label{12713-7}
\end{eqnarray}
\begin{Lemma}\label{theorem1} Suppose Condition 1 holds.  For any  function $f\in\mathcal{D}$ ,  any $i\in \mathcal{S}$, (i)
  the function $R_{f,\pi^{0,b}}(\cdot,i)$ for any  $b\ge 0$, is a continuously differentiable solution on $[0,\infty)$ to the  equations \begin{align}
&\frac{\sigma^2(x,i)}{2}g^{\prime\prime}(x)
+\mu(x,i)g^\prime(x)-(\delta_i+q_i) g(x)+\sum_{j\neq i}q_{ij}f(x,j)=0, \mbox{ $0< x< b,$}\label{1}\\
&\frac{\sigma^2(x,i)}{2}g^{\prime\prime}(x)+
(\mu(x,i)-\bar{l})g^\prime(x)-(\delta_i+q_i) g(x)+\sum_{j\neq i}q_{ij}f(x,j)=-\frac{\bar{l}}{1+d}, \mbox{ $x> b,$}\label{2}\\
&
g^\prime(0+)=\frac1{1-c},\ \ \ \lim_{x\rightarrow\infty}g(x)<\infty,  \label{3}
\end{align} and is twice continuously differentiable on $(0,b)\cup(b,\infty)$;
    (ii)  the function $h_{f,i}(b):=R_{f,\pi^{0,b}}^\prime(b,i)$ is continuous with respect to $b$ for $0<b<\infty$.
\end{Lemma}

Throughout the paper,  we use $\frac{\dif^-}{\dif x}g(x,i)$ and $\frac{\dif^+}{\dif x}g(x,i)$ to represent the derivatives of $g$ from the left- and right-hand side, respectively, with respect to $x$.
 \begin{Corollary}\label{the2} Suppose Condition 1 holds.
  For any $f\in\mathcal{D}$, $i\in\mathcal{S}$ and  $b\ge 0$, (i) $R_{f,\pi^{0,b}}(\cdot,i)$ is  increasing, bounded, continuously differentiable on $(0,\infty)$, and twice continuously differentiable on $(0,b)\cup(b,\infty)$  with $R_{f,\pi^{0,b}}^\prime(0+,i)=\frac1{1-c}$,
 $\left[\frac{\dif^-}{\dif x}R_{f,\pi^{0,b}}^\prime(x,i)\right]_{x=b}=\lim_{x\uparrow b}R_{f,\pi^{0,b}}^{\prime\prime}(x,i)$  and  $ \left[\frac{\dif^+}{\dif x}R_{f,\pi^{0,b}}^\prime(x,i)\right]_{x=b}=\lim_{x\downarrow b}R_{f,\pi^{0,b}}^{\prime\prime}(x,i);$ and  (ii)  if $R_{f,\pi^{0,b}}^\prime(b,i)=\frac1{1+d}$, then $R_{f,\pi^{0,b}}(x,i)$ is twice continuously differentiable with respect to $x$ at $x=b$.
\end{Corollary}

We use $R_{f,\pi^{0,b}}^\prime(0,i)$ and $R_{f,\pi^{0,b}}^{\prime\prime}(0,i)$ to denote $R_{f,\pi^{0,b}}^\prime(0+,i)$ and $R_{f,\pi^{0,b}}^{\prime\prime}(0+,i)$, respectively.

 \begin{Lemma}\label{12713-11} Suppose Conditions 1 and 2 hold. For any fixed $f\in\mathcal{D}$,  $i\in\mathcal{S}$ and  $b\ge 0$, we have $R_{f,\pi^{0,0}}^{\prime\prime}(0+,i)\le 0$, and in the case $b>0$,  $R_{f,\pi^{0,b}}^{\prime\prime}(0+,i)\le 0$ if  $R_{f,\pi^{0,b}}^\prime(b,i)\le \frac1{1-c}$.
 \end{Lemma}

\begin{Lemma}\label{13-2-20} Suppose Conditions 1 and 2 hold.
  For any $f\in\mathcal{D}$ and $i\in\mathcal{S}$, (i) $R_{f,\pi^{0,0}}^{\prime\prime}(x,i)\le 0$ for $x\ge 0$, and in the case  $b> 0$, $R_{f,\pi^{0,b}}^{\prime\prime}(x,i)\le 0$ for $x\ge 0$ if $R_{f,\pi^{0,b}}^\prime(b,i)=\frac1{1+d}$; and
(ii)    for $b>0$, if $R_{f,\pi^{0,b}}^\prime(b,i)>\frac1{1+d}$, $R_{f,\pi^{0,b}}^{\prime\prime}(x,i)\le 0$ for $x\in [ 0,b)$ and $R_{f,\pi^{0,b}}^{\prime\prime }(b-,0)\le 0$.
\end{Lemma}

Let $I\{\cdot\}$ be the indicator function. Define for any fixed $b\ge 0$ and any fixed $\pi\in\Pi$,
\begin{align}
\tau_b^\pi&=\inf\{t\ge 0: X_t^\pi\ge b\},\label{21714-1}\\
{W}_{f,b}(x,i)&=\sup_{\pi\in\Pi}\E_{(x,i)}\Bigg[\int_0^{ \tau_b^\pi\wedge \sigma_1}e^{-\Lambda_s}\frac{l_s}{1+d}\dif s-\int_0^{ \tau_b^\pi\wedge \sigma_1}e^{-\Lambda_s}\frac1{1-c}\dif C_s\nonumber\\
&+e^{-\Lambda_{  \tau_b^\pi}}R_{f,\pi^{0,b}}(X_{\tau_b^\pi}^{\pi},\xi_0)I\{\tau_b^\pi< \sigma_1\}+e^{-\Lambda_{  \sigma_1}}f(X^\pi_{ \sigma_1},\xi_{ \sigma_1})I\{ \sigma_1\le \tau_b^\pi\} \Bigg].\label{450}
\end{align}

\begin{Theorem}\label{the5} Suppose Conditions 1 and 2 hold.
For any $f\in\mathcal{D}$, any $i\in \mathcal{S}$ and any $b>0$,   if $R_{f,\pi^{0,b}}^\prime(b,i)>\frac1{1+d}$,  then  $
  R_{f,\pi^{0,b}}^\prime(x,i)> \frac1{1+d}$  for $0< x\le b$ and   $R_{f,\pi^{0,b}}(x,i)={W}_{f,b}(x,i)$ for $x\ge 0$.
\end{Theorem}

We show in the following theorems  that if $b$ is chosen appropriately, the return function for the strategy $\pi^{0,b}$ coincides with the optimal return function of the modified problem.
\begin{Theorem}\label{23513-9}
Suppose that Conditions 1 and 2 hold. For any  $f\in\mathcal{D}$ and any $i\in \mathcal{S}$, (i) if     $R_{f,\pi^{0,0}}^\prime(0+,i)\le \frac1{1+d}$,  then    $V_f(x,i) =R_{f,\pi^{0,0}}(x,i)$ for $x\ge 0$; and  (ii) if for a fixed  $b>0$, $R_{f,\pi^{0,b}}^\prime(b,i)=\frac1{1+d}$, then $V_f(x,i)=R_{f,\pi^{0,b}}(x,i)$ for $x\ge 0$.
\end{Theorem}

\begin{Lemma}\label{lemff1} Suppose  Conditions 1 and 2 hold,  $f\in\mathcal{D}$ and $i\in \mathcal{S}$. Let $R_{f,\pi^{0,0}}^\prime(0,i)$  denote $R_{f,\pi^{0,0}}^\prime(0+,i)$.
If $R_{f,\pi^{0,b}}^\prime(b,i)>\frac{1}{1+d}$ for all $b\ge 0$, then   $V_f(x,i)=\lim_{b\rightarrow \infty}R_{f,\pi^{0,b}}(x,i)$ for $x\ge 0$.
\end{Lemma}

Again we use $R_{f,\pi^{0,0}}^\prime(0,i)$ to denote $R_{f,\pi^{0,0}}^\prime(0+,i)$.
Define for any $f\in\mathcal{D}$ and $i\in \mathcal{S}$,
 \begin{align}
 b^f_i=\infty\mbox{ if $R_{f,\pi^{0,b}}^\prime(b,i)>\frac1{1+d}$ for all $b\ge 0$, and $b^f_i=
 \inf\{b\ge 0: R_{f,\pi^{0,b}}^\prime(b,i)\le \frac1{1+d}\}$ otherwise.}
 \label{23513-10}
\end{align}
 We show in the following  that the strategy $\pi^{0,b^f_i}$ is optimal for the modified  problem.
.
\begin{Theorem} \label{23513-11} Suppose Conditions 1 and 2 hold. For any $f\in\mathcal{D}$ and any $i\in\mathcal{S}$,
(i)  $0\le b_i^f<\infty$; and
(ii)  $V_f(x,i)=R_{f,\pi^{0,b^f_i}}(x,i)$ for $x\ge 0$.
\end{Theorem}

\section{The Optimality Results}\label{sec4}
We  use the obtained optimality results for the modified optimization problem to address the original optimization problem.  The starting point is to notice that the optimal return function of the  original optimization $V_f$, when the fixed function $f$ is chosen to be the value function of the original optimization, coincides with the value function $V$.
\begin{Theorem} If Conditions 1 and 2 hold,
\label{6214-4}(i) $V\in\mathcal{D}$; (ii) $b_i^V<\infty$ and $V(x,i)=R_{V,\pi^{0,b^V_i}}(x,i)$.
\end{Theorem}
\begin{Theorem} \label{thmj1} Define $\pi^*$ to be the strategy under which,  the dividend pay-out rate at any time $t$ is $\bar{l}I\{X_t^{\pi^*}\}$, and  the company injects capital to maintain the surplus  at level $0$ whenever the surplus tends to go below $0$ without capital injections.  If Conditions 1 and 2 hold, then
  $V(x,i)=V^{\pi^*}(x,i)$ $i\in E$ and  the  strategy $\pi^*$ is an optimal strategy.
\end{Theorem}

\section{Conclusion}
We have addressed the optimal dividend and financing problem for a regime-switching general diffusion model with restricted dividend rates. Our conclusion is that it is optimal to inject capitals only when necessary and at a minimal amount sufficient for the business to continue, and to pay out dividends at the maximal rate, $\bar{l}$, when the  surplus exceeds  the threshold dependent on the environmental state. This result is consistent with the findings for similar problems under simpler model configuration in the literature. For example,  the optimal strategy with restricted dividend rates is of threshold type for the Brownian motion (see \cite{Taksar2000}), the general diffusion (see \cite{Zhu2014b}), and the regime-switching Brownian motion (see \cite{Zhu2014c}).
\numberwithin{equation}{section}
 \renewcommand{\theequation}{A-\arabic{equation}}
  \setcounter{equation}{0}  
  \section*{APPENDIX}  
\subsection*{A.1 \ Proofs for Sections \ref{sec3} and \ref{sec4}}
For any $i\in\mathcal{S}$ and $b\ge0$, define the operator  $\mathcal{B}$ by
\begin{eqnarray}
\mathcal{B}\ g(x,i)=\frac{\sigma^2(x,i)}{2}g^{\prime\prime}(x,i)+
\mu(x,i)g^\prime(x,i)-\delta_i g(x,i).
\end{eqnarray}

\noindent {\bf Proof of  Lemma \ref{complete}}
Consider any   convergent sequence $\{g_n;n=1,2,\cdots\}$ in $\mathcal{D}$ with  limit $g$. It is sufficient to show $g\in\mathcal{D}$. As for any fixed $i$ and $n$, $g_n(\cdot,i)$ is nondecreasing and concave, so is the function $g(\cdot,i)$. The inequality $g(\cdot,i)\le \frac{\bar{l}}{\underline{\delta}(1+d)}$ follows immediately by noticing $g_n(\cdot,i)\le \frac{\bar{l}}{\underline{\delta}(1+d)}$. It remains to show that  $\frac{g(x,i)-g(y,i)}{x-y}\le \frac1{1-c}$ for $0\le x<y$. We use proof by contradiction. Suppose that there exist $x_0$, $y_0$ with $0\le x_0<y_0$  and $j$ such that
$\frac{g(x_0,j)-g(y_0,j)}{x_0-y_0}> \frac1{1-c}$.  Define $\epsilon_0:= \frac12\left( \frac{g(x_0,j)-g(y_0,j)}{x_0-y_0}-\frac1{1-c}\right)$. Clearly,
$\epsilon_0>0$. As $g_n$ converges to $g$,  we can find an $N>0$ such that for all $n\ge N$,
$
||g_n-g||\le \epsilon_0(y_0-x_0).
$
Therefore, $
|g_n(y_0,j)-g(y_0,j)|\le \epsilon_0(y_0-x_0)
$ and $
|g_n(x_0,j)-g(x_0,j)|\le \epsilon_0(y_0-x_0).
$
As a result, $g_n(y_0,j)-g_n(x_0,j)\ge g(y_0,j)- \epsilon_0(y_0-x_0)-(g(x_0,j)+\epsilon_0(y_0-x_0))= g(y_0,j)-g(x_0,j)-2\epsilon_0(y_0-x_0)= \frac{y_0-x_0}{1-c}$.
On the other hand,   we have $\frac{g_n(y_0,j)-g_n(x_0,j)}{y_0-x_0}<\frac{1}{1-c}$ (due to $g_n\in \mathcal{D}$), which is a contradiction.\hfill $\square$

  \noindent {\bf Proof of  Lemma \ref{remff9}}
%
Noting that  $l_s\le \bar{l}$ and that $\sigma_1$ is exponentially distributed with mean $\frac1{q_i}$ and $\Lambda_s=\delta_i s$ for $s\le \sigma_1$, the upper-bounds follow easily from \eqref{jf1}, \eqref{13613-1} and  \eqref{17713-4}.

 Fix $x$ and $y$ with $y>x\ge 0$. Let $\{X_t^x;t\ge 0\}$ and $\{X_t^y;t\ge 0\}$ denote the surplus processes in absence of control with initial surplus $x$ and $y$, respectively. We use  $\pi^x=\{(C^x_t,D^x_t):t\ge 0\}$ with $D^x_t=\int_0^t l_s^x\dif s$ to denote  any admissible control strategy for the process $\{X_t^x;t\ge 0\}$.
Noting that $\{C_t^x;t\ge 0\}$ is right-continuous and increasing, we have the following decomposition:
  $C_t^x=\int_0^t e_s^x \dif s+\sum_{0<s\le t}(C_s^x-C_{s-}^x)$.
 Define $\zeta_0=0$, $\zeta_1=\inf\{s>0:C_{s}^x-C_{s-}^x>0\ \mbox{\ or\ } \xi_s\neq \xi_{s-}\}$ and $\zeta_{n+1}=\{s> \zeta_n:C_{s}^x-C_{s-}^x>0\ \mbox{\ or\ } \xi_s\neq \xi_{s-}\}$ for $n=1,2,\cdots$. Note that $\xi_t=\xi_{\zeta_n}$ for $t\in[\zeta_n,\zeta_{n+1})$ and hence,
 $\dif X_t^{x,
 \pi^x}=(\mu(X_{t-}^{x,\pi^x},\xi_{\zeta_n})-l_t^x+e_t^x)\dif t+\sigma(X_{t-}^{x,\pi^x},\xi_{\zeta_n})\dif W_t$ and $\dif X_t^{y,
 {\pi}^x}=(\mu(X_{t-}^{y,{\pi}^x},\xi_{\zeta_n})
 -{l}_t^x+{e}_t^x)\dif t+\sigma(X_{t-}^{y,{\pi}^x},\xi_{\zeta_n})\dif W_t$  for $t\in(\zeta_n,\zeta_{n+1}), n=0,1,\cdots$.
 By noting $X_0^{x,\pi^x}=X_0^x=x< y=X_0^y=X_0^{y,{\pi}^x}$ and applying the comparison theorem for solutions of stochastic differential equations (see \cite{IkedaWatanabe1977}),  we can show that with probability one, $X_t^{x,\pi^x}\le X_t^{y,
 {\pi}^x}$ for $t\in[0,\zeta_1)$. Further notice that any discontinuity of a surplus process is caused by a jump in the associated process $C^x$ at the same time and hence, $X_{\zeta_1}^{x,\pi^x}
 =X_{\zeta_1-}^{x,\pi^x}+(C_{\zeta_1}^x-C_{\zeta_1-}^x)
 \le X_{\zeta_1-}^{y,\pi^x}+({C}_{\zeta_1}^x-{C}_{\zeta_1-}^x)
 =X_{\zeta_1}^{y,\pi^x}$ with probability one. As a result, by applying the comparison theorem on $(\zeta_1,\zeta_{2})$  we can see  $X_t^{x,\pi^x}\le X_t^{y,
 \pi^x}$ for $t\in(\zeta_1,\zeta_{2})$ with probability one. Repeating the same procedure, we can show that $X_t^{x,\pi^x}\le X_t^{y,
 \pi^x}$  for $t\in (\zeta_n,\zeta_{n+1}]$ with probability one. In conclusion, $X_t^{x,\pi^x}\le X_t^{y,
 \pi^x}$ for all $t\ge 0$  with probability one.  Therefore,  $\pi^x$  satisfies all the requirements for being an admissible strategy for the risk process $X^y$ and hence,
  $R_{f,\pi^x} (y,i)\le  V_f(y,i)$ and $R_{\pi^x} (y,i)\le  V(y,i)$. Using this and \eqref{4214-2} we can show
 $
 R_{f,\pi^x} (x,i)
    \le R_{f,\pi^x} (y,i)\le  V_f(y,i).
 $
Similarly we can obtain
$
 R_{\pi^x} (x,i)
    \le  V(y,i).
 $
 By the arbitrariness of $\pi^x$, we  conclude that $V_f(x,i)\le
     V_f(y,i)$ and  $V(x,i)\le
     V(y,i)$ for $0\le x<y$.
\hfill $\square$

For any  $f\in\mathcal{C}$ and $i\in\mathcal{S}$, define the function
$
w_{f,i}: \mathbb{R}\times\mathcal{S}\rightarrow\mathbb{R}$ by \begin{align}
w_{f,i}(\cdot,i)=R_{f,\pi^{0,b}}(\cdot,i) \mbox{ and } w_{f,i}(\cdot,j)=f(\cdot,j) \mbox{ if $j\neq i$.}
 \end{align}
 \begin{Lemma}\label{itolemma}  For any $f\in\mathcal{C}$ and $i\in\mathcal{S}$, suppose the function  $
w_{f,i}: \mathbb{R}\times\mathcal{S}\rightarrow\mathbb{R}$ with $w_{f,i}(\cdot,j)=f(\cdot,j)$ if $j\neq i$, is bounded, continuously differentiable and piecewise twice continuously differentiable with respect to the first argument on $[0,\infty)$, and the function $w_{f,i}(\cdot,i)$ satisfies the ordinary differential equations \eqref{1} and \eqref{2}. Then, for any  $\pi\in\Pi$,
 there exists a positive sequence of stopping times $\{\tau_n;n=1,2,\cdots\}$ with $\lim_{n\rightarrow \infty}\tau_n=\infty$ such that
  \begin{align}
  &w_{f,i}(x,i)=\E_{(x,i)}\bigg[e^{-\Lambda_{ \tau_n \wedge \sigma_1\wedge t}}w_{f,i}(X^\pi_{ \tau_n\wedge \sigma_1\wedge t}, \xi_{ \tau_n\wedge \sigma_1\wedge t})+\int^{
\tau_n \wedge \sigma_1\wedge t}_0l_se^{-\Lambda_s}w_{f,i}^\prime(X^\pi_{ \tau_n\wedge \sigma_1\wedge t}, \xi_{ \tau_n\wedge \sigma_1\wedge t})\dif s\bigg]
  \nonumber\\
&-\E_{(x,i)}\bigg[\sum_{0<s\le  \tau_n\wedge \sigma_1 \wedge t}e^{-\Lambda_s}\left(w_{f,i}(X^\pi_{s},\xi_{s-})-w_{f,i}(X^\pi_{s-},\xi_{s-})\right)
+\int^{ \tau_n \wedge \sigma_1\wedge t}_0 e^{-\Lambda_s}
w_{f,i}^{\prime}(X^\pi_{s-},\xi_{s-})\dif
\tilde{C}_s\bigg].\nonumber\\
&-\E_{(x,i)}\bigg[\int^{
\tau_n \wedge \sigma_1\wedge t}_0e^{-\Lambda_s}\bar{l}(w_{f,i}^\prime(X^\pi_{s-},\xi_{s-})-\frac{1}{1+d})I\{X^\pi_{s-}\ge b\}\dif s\bigg].
\label{new5}
\end{align}
  \end{Lemma}
  \proof Note that  Applying  It\^o's formula yields that
 \begin{eqnarray}
&&\E_{(x,i)}\bigg[e^{-\Lambda_{ {\tau_n} \wedge \sigma_1\wedge t}}w_{f,i}(X^\pi_{ {\tau_n}\wedge \sigma_1\wedge t}, \xi_{ {\tau_n}\wedge \sigma_1\wedge t}) -
w_{f,i}(X_{0}^{\pi},\xi_0)\bigg]\nonumber\\
&=&I_1+I_2+I_3+\E_{(x,i)}\bigg[\sum_{0<s\le  {\tau_n} \wedge \sigma_1\wedge t}e^{-\Lambda_s}\left(w_{f,i}(X^\pi_{s-},\xi_{s})-w_{f,i}(X^\pi_{s-},\xi_{s-})\right)\bigg],\label{ito}
\end{eqnarray}
where $I_1=\E_{(x,i)}\bigg[\int^{
{\tau_n}\wedge \sigma_1\wedge t}_0e^{-\Lambda_s}\left(\mathcal{B}w_{f,i}(X^\pi_{s-},\xi_{s-})-l_{s}
w_{f,i}^{\prime}(X^\pi_{s-},\xi_{s-}) \right)\dif s \bigg]
$,\\ $I_2=\E_{(x,i)}\bigg[\int^{ {\tau_n}\wedge \sigma_1 \wedge t}_0 e^{-\Lambda_s}\sigma(X^\pi_{s-},\xi_{s-})
w_{f,i}^{\prime}(X^\pi_{s-},\xi_{s-})\dif
W_s\bigg]$ and\\ $I_3=\E_{(x,i)}\bigg[\sum_{0<s\le  {\tau_n}\wedge \sigma_1 \wedge t}e^{-\Lambda_s}\left(w_{f,i}(X^\pi_{s},\xi_{s-})-w_{f,i}(X^\pi_{s-},\xi_{s-})\right)
+\int^{ {\tau_n} \wedge \sigma_1\wedge t}_0 e^{-\Lambda_s}
w_{f,i}^{\prime}(X^\pi_{s-},\xi_{s-})\dif
\tilde{C}_s\bigg].
$

Notice that  the stochastic processes\\ $\int^{t}_0 e^{-\Lambda_s}\sigma(X^{{\pi}}_{s-},\xi_{s-})
w_{f,i}^{\prime}(X^{{\pi}}_{s-},\xi_{s-})\dif
W_s$ and $
\int^{t}_0e^{-\Lambda_s}\left(q_i w_{f,i}(X^{{\pi}}_{s-},\xi_{s-})-\sum_{j\neq i} q_{ij} w_{f,i}(X^{{\pi}}_{s-},j)\right)\dif s+$\\
$\sum_{0<s\le t}e^{-\Lambda_s}\left(
w_{f,i}(X^{{\pi}}_{s-},{{\pi}}_{s})
-w_{f,i}(X^{{\pi}}_{s-},\xi_{s-})\right)$ are $\pr_{(x,i)}$-local martingales. Hence, we can always find a positive sequence of stopping times $\{\tau_n;n=1,2,\cdots\}$ with $\lim_{n\rightarrow \infty}\tau_n=\infty$ such that  both $\int^{t\wedge \tau_n}_0 e^{-\Lambda_s}\sigma(X^{{\pi}}_{s-},\xi_{s-})
w_{f,i}^{\prime}(X^{{\pi}}_{s-},\xi_{s-})\dif
W_s$ and \\
$\int^{ t\wedge \tau_n }_0e^{-\Lambda_s}\left(q_i w_{f,i}(X^{{\pi}}_{s-},\xi_{s-})-\sum_{j\neq i} q_{ij} w_{f,i}(X^{{\pi}}_{s-},j)\right)\dif s$\\
$+\sum_{0<s\le t\wedge \tau_n}
e^{-\Lambda_s}\left(w_{f,i}(X^{{\pi}}_{s-},\xi_{s})
-w_{f,i}(X^{{\pi}}_{s-},\xi_{s-})\right)$ are $\pr_{(x,i)}$-martingales.
Then it follows by the optional stopping theorem that
 \begin{align}
 &I_2=\E_{(x,i)}\bigg[\int^{t\wedge \tau_n\wedge \sigma_1}_0 e^{-\Lambda_s}\sigma(X^{{\pi}}_{s-},\xi_{s-})w_{f,i}^{\prime}
(X^{{\pi}}_{s-},\xi_{s-})\dif
W_s\bigg]=0,\label{new1}\\
&\E_{(x,i)}\bigg[\int^{ t\wedge \tau_n\wedge\sigma_1 }_0e^{-\Lambda_s}\left(q_i w_{f,i}(X^{{\pi}}_{s-},\xi_{s-})-\sum_{j\neq i} q_{ij} w_{f,i}(X^{{\pi}}_{s-},j)\right)\dif s\nonumber\\
 &+\sum_{0<s\le t\wedge \tau_n\sigma_1}
e^{-\Lambda_s}\left(w_{f,i}(X^{{\pi}}_{s-},\xi_{s})
-w_{f,i}(X^{{\pi}}_{s-},\xi_{s-})\right)\bigg]=0.\label{con-1}
\end{align}
Noting that $X_s^{\pi}-X_{s-}^{\pi}=C_s-C_{s-}\ge 0$, $\xi_{s-}=i$, and  $w_{f,i}(X^\pi_{s-}, \xi_{s-})=w_{f,i}(X^\pi_{s-},i)$ for $s\le \sigma_1$ given $\xi_0=i$, that the function $w_{f,i}(\cdot,i) $ satisfies both \eqref{1} and \eqref{2} , and that $w_{f,i}(\cdot,j)=f_i(\cdot,j)$ if $j\neq i$, we obtain that for $s\le {\tau_n}\wedge\sigma_1$,
$
\mathcal{B}w_{f,i}(X^\pi_{s-},\xi_{s-})=q_iw_{f,i}(X^\pi_{s-},\xi_{s-})+
\bar{l}(w_{f,i}^\prime(X^\pi_{s-},\xi_{s-})-\frac{1}{1+d})I\{X^\pi_{s-}\ge b\}-\sum_{j\neq i} q_{ij} w_{f,i}(X^\pi_{s-},j)$, which combined with \eqref{ito},  \eqref{new1}, \eqref{con-1}    and $\E_{(x,i)}\bigg[w_{f,i}(X_{0}^{\pi},\xi_0)\bigg]=w_{f,i}(x,i)$ implies the final result.\hfill $\square$

\noindent {\bf  Proof of Lemma \ref{theorem1}} (i) Let $v_1(\cdot;i)$ and $v_2(\cdot;i)$ denote a set of linearly independent solutions to the equation $
\frac{\sigma^2(x,i)}{2}g^{\prime\prime}(x)+
\mu(x,i)g^\prime(x)-(\delta_i+q_i) g(x)=0,
$ and $v_3(\cdot;i)$ and $v_4(\cdot;i)$  denote       a set of linearly independent solutions  to the equation $
\frac{\sigma^2(x,i)}{2}g^{\prime\prime}(x)+
(\mu(x,i)-\bar{l})g^\prime(x)-(\delta_i+q_i) g(x)=0. $\\
Define
$W_1(x;i)=v_1(x;i)v_2^\prime(x;i)-v_2(x;i)v_1^\prime(x;i)$, $W_2(x;i)=v_3(x;i)v_4^\prime(x;i)-v_4(x;i)v_3^\prime(x;i)$,
$B_1(x;i)=v_1(x;i)\int_0^x\frac{v_2(y;i)}{W_1(y;i)}\frac{2\sum_{j\neq i}q_{ij}f(y,j)}{\sigma^2(y,i)}\dif y-v_2(x;i)\int_0^x\frac{v_1(y;i)}{W_1(y;i)}\frac{2\sum_{j\neq i}q_{ij}f(y,j)}{\sigma^2(y,i)}\dif y$, and
\begin{align*}
&B_2(x;i)=v_3(x;i)\int_0^x\frac{v_4(y;i)}{W_2(y;i)}
\frac{2\left(\bar{l}/(1+d)+\sum_{j\neq i}q_{ij}f(y,j)\right)}{\sigma^2(y,i)}\dif y\nonumber\\
&\ \ \ \ \ \ \ \ \ \ \ \ \ -v_4(x;i)\int_0^x\frac{v_3(y;i)}{W_2(y;i)}\frac{2\left(\bar{l}/(1+d)+\sum_{j\neq i}q_{ij}f(y,j)\right)}{\sigma^2(y,i)}\dif y.
\end{align*}
Then for any constants $K_1,K_2,K_3$ and $K_4$, the functions, $K_1v_1(\cdot;i)+K_2v_2(\cdot;i)+B_1(
\cdot;i),$ and $
K_3v_3(\cdot;i)+K_4v_4(\cdot;i)+B_2(
\cdot;i),
$
are solutions to the equations \eqref{1} and \eqref{2}, respectively. Define the function $g_{b,i}$ by
$
g_{b,i}(x)=
K_1v_1(x;i)+K_2v_2(x;i)+B_1(x;i)$ for $0\le x< b$ and $g_{b,i}(x)=
K_3v_3(x;i)+K_4v_4(x;i)+B_2(x;i)$ for $x\ge b$,
where  $K_1$, $K_2$, $K_3$ and $K_4$ are constants satisfying
\begin{align}
&K_1v_1(b;i)+K_2v_2(b;i)+B_1(b;i)=K_3v_3(b;i)+K_4v_4(b;i)+B_2(b;i),\label{300413-2}\\
&K_1v_1^\prime(b;i)+K_2v_2^\prime(b;i)+B_1^\prime(b;i)=K_3v_3^\prime(b;i)+K_4v_4^\prime(b;i)+B_2^\prime(b;i),\label{300413-03}
\end{align}
\begin{align}
&K_1v_1^\prime(0;i)+K_2v_2^\prime(0;i)=\frac1{1-c},\ \ \lim_{x\rightarrow\infty}(K_3v_3(x;i)+K_4v_4(x;i)+B_2(x;i))<\infty.\label{300413-3}
\end{align}
For  $b\ge 0$, we can easily verify  that $g_{b,i}^\prime(0+)=\frac1{1-c}$, and that $g_{b,i}(\cdot)$ is continuously differentiable on $[0,\infty)$ and twice continuously differentiable on $[0,b)\cup(b,\infty)$. Hence, the existence of a solution with desired property 
 has been proven.

 It suffices to
 show  $R_{f,\pi^{0,b}}(x,i)=g_{b,i}(x)$ for $x\ge 0$. Define $w_{f,i}$  by \begin{eqnarray}
w_{f,i}(x,j)=g_{b,i}(x)\ \mbox{ if $j=i$ and, } w_{f,i}(x,j)=f(x,j) \mbox{ if $j\neq i$}.
\label{25614-100}
\end{eqnarray}
 Note that the process, $X^{0,b}$, will always stay at or above $0$  and the company injects  capital only when the process reaches down to $0$ with a minimal amount to ensure that the surplus never falls below $0$. Further note that $\xi_{s-}=\xi_0$ for $s\le \sigma_1$. Hence, we conclude that the process $C^{0,b}$ is continuous and that given $\xi_0=i$, the following equations hold for $s\le \sigma_1$,
 \begin{eqnarray}
 &&X^{0,b}_{s}=X^{0,b}_{s-}+(C^{0,b}_{s}-C^{0,b}_{s-})=X^{0,b}_{s-},\  \ \ w_i(X^{0,b}_{s},\xi_{s-})-w_i(X^{0,b}_{s-},\xi_{s-})=0\label{17615-4}\\
 &&w_i^{\prime}(X^{0,b}_{s-},\xi_{s-}) \dif \tilde{C}_s^{0,b}=g_{b,i}^{\prime}(0) \dif {C}_s^{0,b}=\frac{\dif {C}_s^{0,b}}{1-c}.\label{23913-1}
\end{eqnarray}
 By  applying  Lemma \ref{itolemma}, we know that for some positive sequence of stopping times $\{\tau_n;n=1,2,\cdots\}$ with $\lim_{n\rightarrow \infty}\tau_n=\infty$, the equation \eqref{new5}  holds. Then  by setting $\pi=\pi^{0,b}$ in \eqref{new5}, noticing that the dividend payment rate at time $s$ is $\bar{l}I\{X^{0,b}_{s-}\ge b\}$ under the strategy $\pi^{0,b}$ and that $g_{b,i}(x)=w_{f,i}(x,i)$, and using  \eqref{17615-4} and  \eqref{23913-1},
we arrive at
\begin{eqnarray}
g_{b,i}(x)&=&\E_{(x,i)}[e^{-\Lambda_{\sigma_1\wedge
\tau_n\wedge t}}w_{f,i}(X^{0,b}_{\sigma_1\wedge\tau_n\wedge t},\xi_{\sigma_1\wedge\tau_n\wedge t})]+\E_{(x,i)}\bigg[\int^{\sigma_1\wedge \tau_n\wedge t}_0 \frac{\bar{l}e^{-\Lambda_s}}{1+d}I\{X^{0,b}_{s}\ge b\}\dif s\bigg]\nonumber\\
&&-\E_{(x,i)}\bigg[\int^{\sigma_1\wedge \tau_n\wedge t}_0\frac{e^{-\Lambda_s}}{1-c} \dif C_s^{0,b}\bigg].\label{12713-6}
\end{eqnarray}
Note that the function $w_{f,i}(
\cdot,\cdot)$ is bounded.
By letting $t\rightarrow\infty$ and $n\rightarrow\infty$ on both sides of \eqref{12713-6}, and then using the dominated convergence for the first expectation on the right-hand side and  the monotone convergence theorem for the other expectations, we can interchange the limits and the expectation and therefore can conclude that
$
g_{b,i}(x)
=R_{f,\pi^{0,b}}(x,i)$ for $x\ge 0.$

\noindent (ii) Note by \eqref{4214-2} that $\lim_{x\rightarrow \infty} g_{b,i}(x)=\lim_{x\rightarrow \infty} R_{f,\pi^{0,b}}(x,i)=
=A_{f,i},
$
 where the second last equality follows by noticing that given $X_0=x$, $X_s^{0,b}\rightarrow \infty$ as $x\rightarrow \infty$ and hence $C_s^{0, b}\rightarrow 0$ as $x\rightarrow \infty$, and the last equality follows by noting that, given $(X_0,\xi_0)=(x,i)$, $\sigma_1$ is exponentially distributed with mean $\frac1{q_i}$, and using the definition of $A_{f,i}$ in \eqref{12713-7}.
 So the constants $K_1,K_2,K_3$ and $K_4$ are solutions to the system of linear equations \eqref{300413-2}-\eqref{300413-3} and $K_3v_3(\infty)+K_4v_4(\infty)+B_2(\infty)=A_{f,i}$. Note that the coefficients of the above system of equations are either constants or continuous functions of $b$. Hence, $K_1,K_2,K_3$ and $K_4$ are continuous functions of $b$, denoted by $K_1(b),K_2(b),K_3(b)$ and $K_4(b)$ here. As a result, the function $h_{f,i}(b)=g_{b,i}^\prime(b)=K_1(b)v_1^\prime(b)+K_2(b)v_2^\prime(b)+B_1^\prime(b;i)$ is continuous  for $0<b<\infty$. \hfill$\square$

For any $f\in\mathcal{C}$, $i\in\mathcal{S}$ and $b\ge 0$, define the functions $h$ and $\bar{h}$ by
\begin{align}
h_{f,i,b}(x)&=(\delta_i+q_i) R_{f,\pi^{0,b}}(x,i)-\mu(x,i)R_{f,\pi^{0,b}}^\prime(x,i)-\sum_{j\neq i}q_{ij}f(x,j)\nonumber\\
&-\bar{l}\left(\frac1{1+d}-R_{f,\pi^{0,b}}^\prime(x,i)\right)I\{x\ge b\},\label{function}
\end{align}
\begin{align}
 \bar{h}_{f,i,b}(x)&=
(\delta_i+q_i) R_{f,\pi^{0,b}}(x,i)-\mu(x,i)R_{f,\pi^{0,b}}^\prime(x,i)-\sum_{j\neq i}q_{ij}f(x,j)\nonumber\\
&-\bar{l}\left(\frac1{1+d}-R_{f,\pi^{0,b}}^\prime(x,i)\right)I\{x> b\}.\label{function1}
 \end{align}

\noindent{\bf  Proof of  Corollary \ref{the2}} (i) is an immediate result of Remark \ref{rem1} and Lemma \ref{theorem1} (i). (ii)  By (i) and Lemma \ref{theorem1}(i) we have  $\left[\frac{\dif^-}{\dif x}R_{f,\pi^{0,b}}^\prime(x,i)\right]_{x=b}
=\lim_{x\downarrow b}\frac{2h_{f,i,b}(b,i)}{\sigma^2(b,i)}
 $ and
 $\left[\frac{\dif^+}{\dif x}R_{f,\pi^{0,b}}^\prime(x,i)\right]_{x=b}
 =\lim_{x\downarrow b}\frac{2h_{f,i,b}(b,i)}{\sigma^2(b,i)}.$ By noting  $R_{f,\pi^{0,b}}^\prime(b,i)=\frac1{1+d},$  we conclude $\left[\frac{\dif^-}{\dif x}R_{f,\pi^{0,b}}^\prime(x,i)\right]_{x=b}=\left[\frac{\dif^+}{\dif x}R_{f,\pi^{0,b}}^\prime(x,i)\right]_{x=b}$. $\square$\\
For any sequence $\{y_n\}$, define
 \begin{align}
 k_{f,b}(x,i;\{y_{n}\})&=
(\delta_i+q_i-\mu^\prime(x,i))R_{f,\pi^{0,b}}^\prime(x,i)-\sum_{j\neq i}q_{ij}\lim_{n\rightarrow \infty} \frac{f(y_{n},j)-f(x,j)}{y_{n}-x}.\label{newdef-1}
\end{align}
\noindent{\bf Proof of    Lemma \ref{12713-11}}
   Throughout the proof, we assume $f\in\mathcal{D}$, $i\in\mathcal{S}$ and $b\ge 0$, unless stated otherwise.  We use proof by contradiction. Suppose $R_{f,\pi^{0,b}}^{\prime\prime}(0+,i)> 0$.

     Since $R_{f,\pi^{0,0}}(\cdot,i)$ is bounded, we can find a large enough $x$ such that $R_{f,\pi^{0,0}}^{\prime}(x,i)<\frac1{1-c}=R_{f,\pi^{0,0}}^{\prime}(0+,i)$, where the last equality is by Lemma \ref{theorem1} (i). Hence there exists an $x>0$ such that $R_{f,\pi^{0,0}}^{\prime\prime}(x,i)< 0$.
 In the case $b>0$, notice that $ R_{f,\pi^{0,b}}^{\prime}(0+,i)=\frac1{1-c}\ge R_{f,\pi^{0,b}}^\prime(b,i).$ So for $b>0$ there exists an $x\in (0,b)$ such that $R_{f,\pi^{0,b}}^{\prime\prime}(x,i)\le 0$.
 Define $x_1=\inf\{x>0:  R_{f,\pi^{0,b}}^{\prime\prime}(x,i)\le 0\}$. Then $x_1>0$ in the case $b=0$ and  $x_1\in(0,b)$ in the case $b>0$, and for $b\ge 0$, \begin{eqnarray}
 R_{f,\pi^{0,b}}^{\prime\prime}(x_1,i)=0,\ \ \mbox{  $R_{f,\pi^{0,b}}^{\prime\prime}(x,i)>0$ for $x\in[0,x_1)$.}\label{12713-9}
   \end{eqnarray}
   As  a result, for $b\ge 0$, \begin{eqnarray}
   R_{f,\pi^{0,b}}^{\prime}(x,i)>R_{f,\pi^{0,b}}^{\prime}(0+,i)=\frac1{1-c}\ \mbox{ for $x\in(0,x_1]$.}\label{12713-10}
    \end{eqnarray} Write $R_{f,\pi^{0,b},i}(x) =R_{f,\pi^{0,b}}(x,i)$. It follows by Lemma \ref{theorem1} that for $b\ge 0$, $\mathcal{A}_{f,i,b}R_{f,\pi^{0,b},i}(x)=0$ for $x>0$.
   %
 Therefore, it follows by \eqref{12713-9} and \eqref{function} that for $b\ge 0$,
 $h_{f,i,b}(x)=
 \frac{\sigma^2(x,i)}{2}R_{f,\pi^{0,b}}^{\prime\prime}(x,i)>0$ for $0<x<x_1$ and $h_{f,i,b}(x_1)
=\frac{\sigma^2(x_1,i)}{2}R_{f,\pi^{0,b}}^{\prime\prime}(x_1,i)=0$.
Hence, we  obtain that for $b\ge 0$,  \begin{eqnarray}
\frac{h_{f,i,b}(x,i)-h_{f,i,b}(x_1,i)}{x-x_1}<0,\ \ 0<x<x_1.\label{12713-12}
 \end{eqnarray}
  Note that $x_1>b$ in the case $b=0$, and  that $x_1<b$ in the case $b>0$. Therefore, we can find a non-negative sequence $\{x_{1n}\}$  with $b<x_{1n}\le x_1$ in the case $b=0$, $x_{1n}\le x_1<b$ in the case $b>0$,  and $\lim_{n\rightarrow \infty}x_{1n}=x_1$ such that $\lim_{n\rightarrow \infty} \frac{f(x_{1n},j)-f(x_1,j)}{x_{1n}-x_1}$ exists.
 By replacing  $x$ in \eqref{12713-12} by $x_{1n}$ and then letting $n\rightarrow \infty$ on both sides of \eqref{12713-12} gives $ k_{f,b}(x_1,i;\{x_{1n}\})-(\mu(x_1,i)-\bar{l}I\{b=0\})R_{f,\pi^{0,b}}^{\prime\prime}(x_1,i)
\ge 0,$
      which combined with \eqref{12713-9} implies
$ \left(\sum_{j\neq i}q_{ij}\lim_{n\rightarrow \infty} \frac{f(x_{1n},j)-f(x_1,j)}{x_{1n}-x_1}
-q_iR_{f,\pi^{0,b}}^\prime(x_1,i)\right)$\\$+\left(\mu^\prime(x_1,i)-\delta_i\right)R_{f,\pi^{0,b}}^\prime(x_1,i)
\le 0.
$
 It follows by this inequality, $R_{f,\pi^{0,b}}^\prime(x_1,i)>\frac1{1-c}$ (see \eqref{12713-10}) and $\lim_{n\rightarrow \infty} \frac{f(x_{1n},j)-f(x_1,j)}{x_{1n}-x_1}\le \frac1{1-c}$ (due to $f\in\mathcal{D}$) that
$\left(\mu^\prime(x_1,i)-\delta_i\right)R_{f,\pi^{0,b}}^\prime(x_1,i)
>0$, which combined with \eqref{12713-10} implies $\mu^\prime(x_1,i)-\delta_i> 0$. This contradicts the assumption that $\mu^\prime(x_1,i)\le \delta_i$ (Condition 2).
 \hfill $\square$

\noindent{\bf   Lemma \ref{13-2-20}}  We consider any fixed $f\in\mathcal{D}$ and $i\in\mathcal{S}$ throughout the proof.
We first show that there exists a positive sequence $\{x_n\}$ with $\lim_{n\rightarrow \infty}x_n=\infty$ such that for $b\ge 0$, \begin{eqnarray}
R_{f,\pi^{0,b}}^{\prime\prime}(x_n,i)\le 0.\label{23513-1}
 \end{eqnarray}
 Suppose the contrary: for some $M>0$, $R_{f,\pi^{0,b}}^{\prime\prime}(x,i)>0$ for all  $x\ge M$. This implies $R_{f,\pi^{0,b}}^\prime(x,i)>R_{f,\pi^{0,b}}^\prime(M+1,i)>R_{f,\pi^{0,b}}^\prime(M,i)\ge 0$ for $x>M+1$, where the last inequality follows by the increasing property of $R_{f,\pi^{0,b}}(\cdot,i)$ (see Corollary \ref{the2}(i)). As a result,
 $
 R_{f,\pi^{0,b}}(x,i)>R_{f,\pi^{0,b}}(M+1,i)+R_{f,\pi^{0,b}}^\prime(M+1,i)(x-M-1)$ for $x>M+1$,
 which
 implies  $\lim_{x\rightarrow \infty} R_{f,\pi^{0,b}}(x,i)=\infty$. This contradicts the boundedness of $R_{f,\pi^{0,b}}(\cdot,i)$ (Corollary \ref{the2}(i)).

Write $R_{f,\pi^{0,b},i}(x)=R_{f,\pi^{0,b}}(x,i).$ By Lemma \ref{theorem1} it follows  that
\begin{eqnarray}
\mathcal{A}_{f,i,b}R_{f,\pi^{0,b},i}(x)=0
\mbox{ for $x> 0$}.\label{13-2-4}
 \end{eqnarray}

\noindent (i) By Lemma  \ref{theorem1} and Corollary \ref{the2} we can see that
$R_{f,\pi^{0,b},i}(\cdot)$ is  twice continuously differentiable on $[0,\infty)$ with the differentiability at $0$ referring to the differentiability from the right-hand side.
 It follows by noting $R_{f,\pi^{0,b},i}^{\prime}(b)=R_{f,\pi^{0,b}}^{\prime}(b,i)=\frac1{1+d}\le \frac1{1-c}$ for $b>0$, and   Lemma \ref{12713-11} that
\begin{eqnarray}
R_{f,\pi^{0,b},i}^{\prime\prime}(0+)\le 0\ \mbox{ for $b\ge 0$.} \label{13-2-1}
\end{eqnarray}

We use proof by contradiction to prove the statement in (i).
Suppose that  the statement in (i)  is not true. Then there  exists  a  $b\ge 0$ and a $y_0>                                                                                                                                                              0$ such that $R_{f,\pi^{0,b},i}^{\prime\prime}(y_0)=R_{f,\pi^{0,b}}^{\prime\prime}(y_0,i)>0$. Let $\{x_n\}$ be the sequence defined as before. We can find a positive integer  $N$ such that $x_N>y_0$.  By noting $R_{f,\pi^{0,b},i}^{\prime\prime}(x_N)=R_{f,\pi^{0,b}}^{\prime\prime}(x_N,i)\le 0$ (due to \eqref{23513-1}), \eqref{13-2-1} and the continuity of $R_{f,\pi^{0,b},i}^{\prime\prime}(\cdot)$, we can find $y_1, y_2$ with  $0\le y_1<y_0<y_2\le x_N$ such that \begin{eqnarray}
R_{f,\pi^{0,b}}^{\prime\prime}(y_1,i)=0,\ \ R_{f,\pi^{0,b}}^{\prime\prime}(y_2,i)=0,\ \ \mbox{and }\ R_{f,\pi^{0,b}}^{\prime\prime}(x,i)>0\ \mbox{ for $x\in (y_1,y_2)$.} \label{13-2-5}
 \end{eqnarray}
 Hence,
\begin{eqnarray}
R_{f,\pi^{0,b},i}^{\prime}(y_2)> R_{f,\pi^{0,b},i}^{\prime}(y_1).\label{13-2-2}
\end{eqnarray}
It follows by \eqref{13-2-4} and \eqref{function} that  $-\frac{\sigma^2(x,i)}2R_{f,\pi^{0,b},i}^{\prime\prime}(x)=h_{f,b,i}(x)$ for $x> 0$.
 Note that for $x>0$, $I\{x\ge b\}=I\{x>b\}$ in the case $b=0$, and that in the case $b>0$,  $\frac1{1+d}-R_{f,\pi^{0,b}}^\prime(b,i)=0$ and hence, $\bar{l}\left(\frac1{1+d}-R_{f,\pi^{0,b}}^\prime(x,i)\right)I\{x\ge b\}=\bar{l}\left(\frac1{1+d}-R_{f,\pi^{0,b}}^\prime(x,i)\right)I\{x> b\}$ for $x>0$.
 Therefore, for $x> 0$,
$
\frac{\sigma^2(x,i)}2R_{f,\pi^{0,b}}^{\prime\prime}(x,i)=\bar{h}_{f,i,b}(x)$, which combined with \eqref{13-2-5} implies that for $x\in(y_1,y_2)$,
\begin{eqnarray}
&&\bar{h}_{f,i,b}(y_1)
=\frac{\sigma^2(y_1,i)}2R_{f,\pi^{0,b}}^{\prime\prime}(y_1,i)
=0
<\frac{\sigma^2(x,i)}2R_{f,\pi^{0,b}}^{\prime\prime}(x,i)=
\bar{h}_{f,i,b}(x),\label{13-2-8}\\
&&\bar{h}_{f,i,b}(y_2)
=\frac{\sigma^2(y_2,i)}2R_{f,\pi^{0,b}}^{\prime\prime}(y_2,i)
=0
<\frac{\sigma^2(x,i)}2R_{f,\pi^{0,b}}^{\prime\prime}(x,i)=
\bar{h}_{f,i,b}(x). \label{13-2-9}
 \end{eqnarray}
 Let $\{y_{1n}\}$ and $\{y_{2n}\}$ be two sequences with $y_{1n}\downarrow y_1$ and $y_{2n}\uparrow y_2$ as $n\rightarrow \infty$ such that $\lim_{n\rightarrow\infty}
 \frac{f(y_{1n},j)-f(y_1,j)}{y_{1n}-y_1}$  and $\lim_{n\rightarrow\infty}\frac{f(y_{2n},j)-f(y_2,j)}
 {y_{2n}-y_2}$ exist for all $j\in \mathcal{S}$. It follows by \eqref{13-2-8} and \eqref{13-2-9} that $\frac{\bar{h}_{f,i,b}(y_{1n})
 -\bar{h}_{f,i,b}(y_{1})}{y_{1n}-y_1}>0>
  \frac{\bar{h}_{f,i,b}(y_{2n})-\bar{h}_{f,i,b}(y_{2})}{y_{2n}-y_2}
$.
 By  letting $n\rightarrow \infty$,  we obtain\\
 $k_{f,b}(y_1,i;\{y_{1n}\})-\mu(y_1,i)R_{f,\pi^{0,b}}^{\prime\prime}(y_1,i)
  +\bar{l}R_{f,\pi^{0,b}}^{\prime\prime}(y_1,i)I\{y_1> b\}\ge 0$\\ and
  $k_{f,b}(y_2,i;\{y_{2n}\})
  -\mu(y_2,i)R_{f,\pi^{0,b}}^{\prime\prime}(y_2,i)
  +\bar{l}R_{f,\pi^{0,b}}^{\prime\prime}(y_2,i)I\{y_2> b\} \le 0. $
Therefore, by noting $R_{f,\pi^{0,b}}^{\prime\prime}(y_1,i)=0=R_{f,\pi^{0,b}}^{\prime\prime}(y_2,i)$ (see \eqref{13-2-5}) we have
\begin{eqnarray}
&&k_{f,b}(y_1,i;\{y_{1n}\})\ge 0\ge k_{f,b}(y_2,i;\{y_{2n}\}).
 \label{13-2-14}
 \end{eqnarray}
 On the other hand, note that $0<\delta_i+q_i-\mu^\prime(y_1,i)\le\delta_i+q_i-\mu^\prime(y_2,i)$ (due to the concavity of $\mu(\cdot,i)$), $R_{f,\pi^{0,b}}^{\prime}(y_1,i)<R_{f,\pi^{0,b}}^{\prime}(y_2,i)$ (see \eqref{13-2-2}), $\lim_{n\rightarrow\infty}\frac{f(y_{1n},j)-f(y_1,j)}{y_{1n}-y_1}\ge\lim_{n\rightarrow\infty}\frac{f(y_{2n},j)-f(y_2,j)}{y_{2n}-y_2}$ (due to the concavity of $f(\cdot,j)$). As a result,
 $
  k_{f,b}(y_1,i;\{y_{1n}\})< k_{f,b}(y_2,i;\{y_{2n}\})$, which is a contradiction to \eqref{13-2-14}.

\noindent (ii) We distinguish two cases: (a) $R_{f,\pi^{0,b}}^{\prime\prime}(b+,i)>0$ and (b) $R_{f,\pi^{0,b}}^{\prime\prime}(b+,i)\le  0$.\\
(a) Suppose $R_{f,\pi^{0,b}}^{\prime\prime}(b+,i)>0$. By \eqref{23513-1} we can find $N>0$ such that $x_N>b$ and $R_{f,\pi^{0,b}}^{\prime\prime}(x_N,i)\le 0$. Then by the continuity of the function $R_{f,\pi^{0,b}}^{\prime\prime}(\cdot,i)$ on $(b,\infty)$ (see Corollary \ref{the2}(i))   we know that there exists a $y_2\in (b,x_N]$ such that
   $
   R_{f,\pi^{0,b}}^{\prime\prime}(y_2,i)=0\ \mbox{ and }\ R_{f,\pi^{0,b}}^{\prime\prime}(x,i)>0\ \mbox{ for $x\in(b,y_2)$}.$
   We now proceed to show that $R_{f,\pi^{0,b}}^{\prime\prime}(b-,i)\le 0$. Suppose the contrary, i.e., $R_{f,\pi^{0,b}}^{\prime\prime}(b-,i)> 0$. By noting $R_{f,\pi^{0,b}}^{\prime\prime}(0+,i)\le 0$ (see \eqref{13-2-1}), it follows that there exists a $y_1\in (0, b)$ such that
   $
   R_{f,\pi^{0,b}}^{\prime\prime}(y_1,i)=0$  and $R_{f,\pi^{0,b}}^{\prime\prime}(x,i)>0$  for $x\in(y_1,b)$. In summary, \eqref{13-2-5} holds for $x\in(y_1,y_2)-\{b\}$.
   Rrepeating the argument right below \eqref{13-2-5} in (i), we obtain a contradiction.

     \noindent (b)  Suppose $R_{f,\pi^{0,b}}^{\prime\prime}(b+,i)\le 0$. It follows by \eqref{13-2-4} and the assumption $R_{f,\pi^{0,b}}^\prime(b,i)>\frac1{1+d}$ that
  \begin{eqnarray}
&&R_{f,\pi^{0,b}}^{\prime\prime}(b-,i)=\lim_{x\uparrow b}\frac{2h_{f,i,b}(x,i)}{\sigma^2(x,i)}< \lim_{x\downarrow b}\frac{2h_{f,i,b}(x,i)}{\sigma^2(x,i)}
=R_{f,\pi^{0,b}}^{\prime\prime}(b+,i)\le 0.\label{23513-7}
 \end{eqnarray}
  We now  show that $R_{f,\pi^{0,b}}^{\prime\prime}(x,i)\le 0$ for all $x\in[0,b)$. Suppose the contrary. That is, there exists some $x\in[0,b)$ such that $R_{f,\pi^{0,b}}^{\prime\prime}(x,i)> 0$. Then by noting $R_{f,\pi^{0,b}}^{\prime\prime}(0+,i)\le 0$ (see \eqref{13-2-1}) and $R_{f,\pi^{0,b}}^{\prime\prime}(b-,i)< 0$ (see \eqref{23513-7}), we can find $y_1$ and $y_2$ with $0\le y_1<y_2< b$ such that
  $
   R_{f,\pi^{0,b}}^{\prime\prime}(y_1,i)=0$, $R_{f,\pi^{0,b}}^{\prime\prime}(y_2,i)=0$ and  $R_{f,\pi^{0,b}}^{\prime\prime}(x,i)>0$  for $x\in(y_1,y_2)$.
      Repeating again the argument  right after \eqref{13-2-5}  in (i), we can obtain a contradiction.
     \hfill $\square$

\noindent{\bf   Theorem \ref{the5}}  Note that $\tau_b^\pi=0$ given $X_0^\pi\ge b$. Hence, it follows  from the definition \eqref{450} that \begin{eqnarray}
{W}_{f,b}(x,i)=\sup_{\pi\in\Pi}\E_{(x,i)}\left[R_{f,\pi^{0,b}}(X^\pi_{0},\xi_0) \right]=R_{f,\pi^{0,b}}(x,i)\ \mbox{ for $x\ge b$ and $b=0$}.\label{13-2-40}
\end{eqnarray}
 We consider the case $b>0$. By Lemma \ref{13-2-20} (ii) we know that $R_{f,\pi^{0,b}}^{\prime\prime}(x,i)\le 0$ for $x\in [0,b)$, and $R_{f,\pi^{0,b}}^{\prime\prime}(b-,i)\le 0$.
  Therefore, it follows by  Corollary \ref{the2}(i) that
  \begin{eqnarray}
  \frac1{1-c}= R_{f,\pi^{0,b}}^\prime(0+,i)\ge R_{f,\pi^{0,b}}^\prime(x,i)\ge R_{f,\pi^{0,b}}^\prime(b,i)> \frac1{1+d}\ \mbox{ for $0<x\le b$}.\label{42}
   \end{eqnarray}
 Define $
 w_{f,i}(y,j)=
 R_{f,\pi^{0,b}}(y,i) \mbox{ if $j=i$, and $w_{f,i}(y,j)=
 f(y,j)$ if $j\neq i$}.
$ Then by Corollary \ref{the2}(i) and Lemma \ref{theorem1} we know that $w_i(\cdot,j)$ satisfies the conditions in Lemma \ref{itolemma}.  Then by applying Lemma \ref{itolemma}  we know that for some positive sequence of stopping times $\{\tau_n;n=1,2,\cdots\}$ with $\lim_{n\rightarrow \infty}\tau_n=\infty$, the equation \eqref{new5} holds. By letting $t$ in \eqref{new5} be $\tau_b^\pi\wedge t$,  noting that $X_s^{\pi}-X_{s-}^{\pi}=C_s-C_{s-}\ge 0$, and that given $(X_0,\xi_0)=(x,i)$,  $X^\pi_{s-}\in[0,b)$ and $w_i(X^\pi_{s-}, \xi_{s-})=R_{f,\pi^{0,b}}(X^\pi_{s-},i)$ for $s\le \sigma_1\wedge \tau_b^\pi$, that $\sum_{0<s\le  \tau_n\wedge \sigma_1 \wedge \tau_b^\pi \wedge t}e^{-\Lambda_s}\frac{X^\pi_{s}-X^\pi_{s-}}{1-c}
+\int^{\tau_n\wedge  \sigma_1 \tau_b^\pi \wedge\wedge t}_0
\frac{e^{-\Lambda_s}}{1-c}\dif
\tilde{C}_s=
\int_{0}^{\tau_n\wedge  \sigma_1 \tau_b^\pi \wedge\wedge t}\frac{e^{-\Lambda_s}}{1-c}\dif
C_s$, and using\eqref{42}, we derive that for any $\pi\in \Pi$,  $t>0$ and $0\le x\le b$,
 \begin{eqnarray}
&&\E_{(x,i)}\bigg[\int_0^{\tau_n\wedge  \sigma_1\wedge \tau_b^\pi \wedge t}\frac{l_se^{-\Lambda_s}}{1+d}\dif s-\int_0^{\tau_n\wedge  \sigma_1\wedge \tau_b^\pi \wedge t}\frac{e^{-\Lambda_s}}{1-c}\dif C_s\nonumber\\&+&e^{-\Lambda_{\tau_n\wedge  \sigma_1\wedge \tau_b^\pi \wedge t}}w_i(X^\pi_{\tau_n\wedge  \sigma_1\wedge \tau_b^\pi \wedge t},\xi_{\tau_n\wedge  \sigma_1\wedge \tau_b^\pi \wedge t}) \Bigg]
\le
R_{f,\pi^{0,b}}(x,i).\label{45001}
\end{eqnarray}
 Note that the functions  $R_{f,\pi^{0,b}}(\cdot,j)$ and $f(\cdot,j)$ $j\in\mathcal{S}$ are all bounded. Hence, the functions $w_i(\cdot,j)$ $j\in\mathcal{S}$ are also bounded.
 By letting $\tau_n\rightarrow\infty$ and $t\rightarrow \infty$ on both sides of \eqref{45001},   using the monotone convergence theorem and the dominated convergence theorem and noticing   that due to   $\xi_{s}=\xi_0$ for $0\le s<\sigma_1$ we have
 $\E_{(x,i)}\bigg[e^{-\Lambda_{  \tau_b^\pi\wedge\sigma_1}}w_{f,i}(X^\pi_{ \tau_b^\pi\wedge \sigma_1},\xi_{ \tau_b^\pi\wedge \sigma_1})\bigg]=\E_{(x,i)}\bigg[e^{-\Lambda_{ \tau_b^\pi}}R_{f,\pi^{0,b}}(b,\xi_0)I\{\tau_b^\pi< \sigma_1\}+e^{-\Lambda_{ \sigma_1}}f(X_{\sigma_1},\xi_{\sigma_1})I\{\sigma_1\le \tau_b^\pi\}\bigg]$ and that $\pi$ is an arbitrary admissible strategy and  \eqref{450},
 we can conclude
\begin{eqnarray}
W_{f,b}(x,i)\le
R_{f,\pi^{0,b}}(x,i)\ \mbox{ for } 0\le x\le b.\label{4220}
\end{eqnarray}
Note that   $\{(X_t^{0,b},\xi_t);t\ge 0\}$ is a strong Markov process and that by the Markov property  it follows that
\begin{align}
R_{f,\pi^{0,b}}(x,i)
&=\E_{(x,i)}\bigg[\int_0^{  \tau_b^{\pi^{0,b}}\wedge \sigma_1}\frac{\bar{l}e^{-\Lambda_s}}{1+d}I\{X^{0,b}_s\ge b\}\dif s-\int_0^{  \tau_b^{\pi^{0,b}}\wedge \sigma_1}\frac{e^{-\Lambda_s}}{1-c}\dif C_s\nonumber\\
&+e^{-\delta( \tau_b^{\pi^{0,b}}\wedge \sigma_1)}R_{f,\pi^{0,b}}(X_{  \tau_b^{\pi^{0,b}}\wedge \sigma_1}^{0,b},\xi_{  \tau_b^{\pi^{0,b}}\wedge \sigma_1}) \bigg]\le W_{f,b}(x,i) \mbox{  for $ x\ge 0$},\label{510}
\end{align}
where the last inequality follows by noting $\pi^{0,b}\in \Pi$ and the definition \eqref{450}.

 Combining \eqref{13-2-40},  \eqref{4220} and \eqref{510} completes the proof.
\hfill$\square$

\noindent {\bf Proof of  Theorem \ref{23513-9}}
  We first show that
\begin{eqnarray}
 && R_{f,\pi^{0,b}}^\prime(x,i)\le R_{f,\pi^{0,b}}^\prime(b,i)= \frac1{1+d}\ \mbox{ for $x>b$, $b\ge 0$}.\label{13-2-31}
\end{eqnarray}By Lemma \ref{13-2-20}(i) it follows that $R_{f,\pi^{0,0}}^{\prime\prime}(x,i)\le 0$ for $x\ge 0$. As a result, \eqref{13-2-31} holds for $b=0$. Now suppose $b>0$. By Lemma \ref{theorem1} (i) we know that $R_{f,\pi^{0,b}}^\prime(0+,i)=\frac1{1-c}$. Since $R_{f,\pi^{0,b}}^\prime(b,i)=\frac1{1+d}$, it follows by Corollary \ref{the2} (ii) that $R_{f,\pi^{0,b}}(\cdot,i)$ is twice continuously differentiable on $[0,\infty)$ and  by Lemma \ref{13-2-20} (i) that $R^{\prime\prime}_{f,0,b}(x,i)\le 0$ for $x\ge 0$. Hence, \eqref{13-2-31} holds for $b>0$ as well, and
\begin{eqnarray}
&&\frac1{1-c}=R_{f,\pi^{0,b}}^\prime(0+,i)\ge R_{f,\pi^{0,b}}^\prime(x,i)\ge R_{f,\pi^{0,b}}^\prime(b,i)= \frac1{1+d}\ \mbox{ for $x\in [ 0,b]$}.\label{17713-5}
\end{eqnarray}

It follows by using \eqref{13-2-31} and \eqref{17713-5}, and noting $\bar{l}\ge l_s$ for $s\ge 0$ we obtain that for $b\ge 0$,
\begin{eqnarray}
&&\bar{l}I\{X^\pi_{s}\ge b\}\left(R_{f,\pi^{0,b}}^{\prime}(X^\pi_{s-},i)-\frac1{1+d}\right)
-l_sR_{f,\pi^{0,b}}^{\prime}(X^\pi_{s-},i)\nonumber\\
&=&(\bar{l}-l_s)I\{X^\pi_{s}\ge b\}R_{f,\pi^{0,b}}^{\prime}(X^\pi_{s-},i)-\frac{\bar{l}}{1+d}I\{X^\pi_{s}\ge b\}-l_sI\{X^\pi_{s}< b\}R_{f,\pi^{0,b}}^{\prime}(X^\pi_{s-},i)\nonumber\\
&\le&\frac{\bar{l}-l_s}{1+d}I\{X^\pi_{s}\ge b\}-\frac{\bar{l}}{1+d}I\{X^\pi_{s}\ge b\}-\frac{l_s}{1+d}I\{X^\pi_{s}< b\}=-\frac{l_s}{1+d},\label{0fs9}
\end{eqnarray}
By \eqref{13-2-31} again we can obtain
\begin{eqnarray}
 && R_{f,\pi^{0,b}}^\prime(x,i)\le \frac1{1-c}\ \mbox{ for $b\ge 0$ and $x>b$}.\label{13-2-313}
\end{eqnarray}
Further, note that for $b\ge 0$ and any  $t\ge 0$,
\begin{align}
&\E_{(x,i)}\bigg[\int_{0<s\le  \sigma_1 \wedge
t}e^{-\Lambda_s}R_{f,\pi^{0,b}}^\prime(X^\pi_{s},\xi_{s-})\dif \tilde{C}_s
+ \sum_{0<s\le  \sigma_1 \wedge
t}e^{-\Lambda_s}\left(R_{f,\pi^{0,b}}(X^\pi_{s},\xi_{s-})-R_{f,\pi^{0,b}}(X^\pi_{s-},\xi_{s-})\right)\bigg]\nonumber\\
&\le \E_{(x,i)}\bigg[\int_{0}^{ \sigma_1 \wedge
t}\frac{e^{-\Lambda_s}}{1-c}\dif \tilde{C}_s+ \sum_{0<s\le  \sigma_1 \wedge
t}\frac{e^{-\Lambda_s}}{1-c}(X^\pi_{s}-X^\pi_{s-})\bigg]
=\E_{(x,i)}\bigg[ \sum_{0<s\le \sigma_1 \wedge
t}\frac{e^{-\Lambda_s}}{1-c}\dif C_s\bigg],\label{17713-7}
\end{align}
where the last inequality follows by \eqref{17713-5}, \eqref{13-2-313}, $\dif\tilde{C}_s\ge 0$, $X_s^{\pi}-X_{s-}^{\pi}=C_s-C_{s-}\ge 0$ and $\dif C_s=\dif \tilde{C}_s+C_s-C_{s-}$.

Define $
w_{f,i}(y,j)=R_{f,\pi^{0,b}}(y,i)$ if $j=i$, and $w_{f,i}(y,j)=f(y,j)$  if $j\neq i$.
Then by Corollary \ref{the2}(i) and Lemma \ref{theorem1} we know that the conditions in Lemma \ref{theorem1} are satisfied. By applying Lemma \ref{itolemma}  we know that for some positive sequence of stopping times $\{\tau_n;n=1,2,\cdots\}$ with $\lim_{n\rightarrow \infty}\tau_n=\infty$, the equation \eqref{new5} holds for any $\pi\in\Pi$, any $b,t>0$ and any $n\in\mathbb{N}$. By using
 \eqref{new5}, \eqref{0fs9}  and \eqref{17713-7} (setting $t=t\wedge \tau_n$) we arrive at
$R_{f,\pi^{0,b}}(x,i)\ge \E_{(x,i)}\bigg[\int^{ \sigma_1\wedge
t\wedge \tau_n}_0\frac{l_se^{-\Lambda_s}}{1+d}\dif s-\sum_{0}^{ \sigma_1 \wedge
t\wedge \tau_n}\frac{e^{-\Lambda_s}}{1-c}\dif C_s
+e^{-\Lambda_{ \sigma_1\wedge
t\wedge \tau_n}}w_{f,i}(X^\pi_{ \sigma_1\wedge t\wedge \tau_n},\xi_{ \sigma_1\wedge t\wedge \tau_n})\bigg]$ for $b\ge 0$.
By noting  that the functions $R_{f,\pi^{0,b}}(\cdot,i)$ and $f(\cdot,j)$, $j\in \mathcal{S}$ are  bounded and letting $t\rightarrow \infty$ and then $n\rightarrow \infty$ and then using the monotone convergence theorem for the first two terms inside the expectation  and the dominated convergence theorem for the last term, we  obtain that for $b\ge 0$,
$
R_{f,\pi^{0,b}}(x,i)
\ge\E_{(x,i)}\bigg[\int^{ \sigma_1}_0\frac{l_se^{-\Lambda_s}}{1+d}\dif s-\int_0^{ \sigma_1 }\frac{e^{-\Lambda_s}}{1-c}\dif C_s
+e^{-\Lambda_{ \sigma_1}}w_{f,i}(X^\pi_{ \sigma_1},\xi_{ \sigma_1})\bigg]
.$
By noting $w_{f,i}(X^\pi_{ \sigma_1},\xi_{ \sigma_1})=f(X^\pi_{ \sigma_1},\xi_{ \sigma_1})$ given $\xi_0=i$, the arbitrariness of $\pi$ and the definition of $V_f$ in \eqref{17713-4} we conclude
$
R_{f,\pi^{0,b}}(x,i)\ge
V_f(x,i)\ \mbox{ for $x\ge 0$.}$
On the other hand, $R_{f,\pi^{0,b}}(x,i)\le V_f(x,i)$ for $x\ge 0$ according to the definition \eqref{17713-4}. Consequently, $R_{f,\pi^{0,b}}(x,i)=V_f(x,i)$ for $x\ge 0$.
\hfill$\square$

\noindent {\bf Proof of  Lemma \ref{lemff1}}  Recall that $\tau_b^\pi$ is defined in \eqref{21714-1}. By Theorem \ref{the5} it follows that for any large enough $b$ and any $x\ge 0$,
\begin{align*}
R_{f,\pi^{0,b}}(x,i)&=W_{f,b}(x,i)=\sup_{\pi\in\Pi}\E_{(x,i)}\bigg[\int^{ \sigma_1 \wedge
\tau_b^\pi}_0\frac{l_se^{-\Lambda_s}}{1+d}\dif s-\int^{ \sigma_1 \wedge
\tau_b^\pi}_0\frac{e^{-\Lambda_s}}{1-c}\dif C_s\nonumber\\
&+e^{-\Lambda_{\tau_b^\pi}}R_{f,\pi^{0,b}}(b, \xi_0)I\{\tau_b^\pi< \sigma_1\}+e^{-\Lambda_{\sigma_1}}f(X^\pi_{\sigma_1}, \xi_{\sigma_1})I\{\sigma_1\le \tau_b^\pi\} \bigg] \nonumber\\
\ge& \sup_{\pi\in\Pi}\E_x\Bigg[\int^{ \sigma_1\wedge
\tau_b^\pi}_0\frac{l_se^{-\Lambda_s}}{1+d}\dif s-\int^{ \sigma_1 \wedge
\tau_b^\pi}_0\frac{e^{-\Lambda_s}}{1-c}\dif C_s+e^{-\Lambda_{\sigma_1}}f(X^\pi_{\sigma_1}, \xi_{\sigma_1})I\{\sigma_1\le \tau_b^\pi \}\Bigg].
\end{align*}
Note  $\lim_{b\rightarrow\infty}\tau_b^\pi=\infty$ and $f$ is bounded. Then it follows by letting $b\rightarrow \infty$ on both sides, and then using the monotone convergence theorem twice and the dominated convergence  that
$\liminf_{b\rightarrow\infty}R_{f,\pi^{0,b}}(x,i)\ge\sup_{\pi\in\Pi}\E_{(x,i)}\left[\int^{ \sigma_1 }_0\frac{l_se^{-\Lambda_s}}{1+d}\dif s-\int^{ \sigma_1 }_0\frac{e^{-\Lambda_s}}{1-c}\dif C_s+e^{-\Lambda_{\sigma_1}}f(X^\pi_{\sigma_1}, \xi_{\sigma_1})\right]=V_f(x,i)$  for $x\ge 0.$
This combined with the fact $R_{f,\pi^{0,b}}(x,i)\le V_f(x,i)$ for $x\ge 0$ completes the proof.\hfill $\square$

\noindent{\bf Proof of  Theorem \ref{23513-11}}
  (i)  $b_i^f\ge 0$  is obvious by the definition. We just need to prove $b_i^f<\infty$. Suppose the contrary. Then by \eqref{23513-10} we have $R_{f,\pi^{0,b}}^{\prime}(b,i)> \frac1{1+d}$ for all $b\ge 0$. Hence, it follows by Lemma \ref{lemff1} that $V_f(x,i)=\lim_{b\rightarrow \infty}R_{f,\pi^{0,b}}(x,i)$ for $x\ge 0$.
For any $b\ge 0$, by Theorem \ref{the5} we know $R_{f,\pi^{0,b}}^\prime(x,i)>\frac1{1+d}$ for $x\in(0,b]$, which implies
$R_{f,\pi^{0,b}}(x,i)>R_{f,\pi^{0,b}}(0,i)+\frac{x}{1+d}$ for $x\in(0,b]$. Hence, for any $x\ge 0$, we can find a $b>x$ such that $V_f(x,i)\ge R_{f,\pi^{0,b}}(x,i)>R_{f,\pi^{0,b}}(0,i)+\frac{x}{1+d}$. Hence, $\lim_{x\rightarrow\infty}V_f(x,i)=+\infty$, which contradicts $V_f(x,i)\le \frac{\bar{l}}{\underline{\delta}(1+d)}$ for $x\ge 0$ (see Lemma \ref{remff9}).
\noindent (ii) is a result of (i) and Theorem \ref{23513-9}.
\hfill $\square$

\noindent {\bf Proof of  Theorem \ref{6214-4}}
(i) Define an operator $\mathcal{P}$ by
\begin{eqnarray}
\mathcal{P}(f)(x,i):=V_f(x,i),\ x\ge 0,i\in\mathcal{S}\ \mbox{ and $f\in\mathcal{C}$}.\label{6214-2}
\end{eqnarray}
Then by Theorem \ref{23513-11} we have, \begin{eqnarray}
\mathcal{P}(f)(x,i)=V_f(x,i)=R_{f,\pi^{0,b_i^f}}(x,i),\ x\ge 0,i\in\mathcal{S}\ \mbox{ and $f\in\mathcal{C}$}.\label{6214-1}
\end{eqnarray}

Recall that $\mathcal{D}\subset\mathcal{C}$ and $(\mathcal{D},||\cdot||)$ is a complete space.
 We will first show that $\mathcal{P}$ is a contraction on  $(\mathcal{D},||\cdot||)$. Consider any $f\in\mathcal{D}$. It follows by Lemma \ref{remff9} and \eqref{6214-1} that  $\mathcal{P}(f)=V_f\in\mathcal{C}$. Note that for any $f\in\mathcal{D}$ and $i\in \mathcal{S}$, $b_i^f<\infty$ according to Theorem \ref{23513-11}. Further notice that by Lemma \ref{theorem1} (ii), we know $R_{f,\pi^{0,b}}^\prime(b,i)$ is continuous in $b$ and $R_{f,\pi^{0,0}}^\prime(0+,i)=\frac1{1-c}>\frac1{1+d}$ by Corollary \ref{the2} (i).  Hence, according to the definition of $b_i^f$ in \eqref{23513-10}, we  have  $R_{f,\pi^{0,b^f_i}}^\prime(b_i^f,i)=\frac1{1+d}$.  Therefore, it follows by Corollary \ref{the2} that
for any $i\in\mathcal{ S}$, the function $R_{f,\pi^{0,b^f_i}}(\cdot,i)$ is twice continuously differentiable on $(0,\infty)$ and by Lemma \ref{13-2-20} (i) that $R_{f,\pi^{0,b^f_i}}(\cdot,i)$ is concave. Notice that by Corollary \ref{the2} (i) again $R_{f,\pi^{0,b^f_i}}^\prime(0+,i)=\frac1{1-c}$. Hence, $\frac{\dif}{\dif x}\mathcal{P}(f)(x,i)=R_{f,\pi^{0,b^f_i}}^\prime(x,i)\le R_{f,\pi^{0,b^f_i}}^\prime(0+,i)=\frac1{1-c}$ for $x>0$, which results in $\frac{\mathcal{P}(f)(x,i)-\mathcal{P}(f)(y,i)}{x-y}\le \frac1{1-c}$ for $0\le x<y$. Therefore, we can conclude $\mathcal{P}(f)\in\mathcal{D}$. For any $f_1,f_2\in\mathcal{D}$, it follows by \eqref{6214-2}  that
\begin{eqnarray}
 && ||\mathcal{P}(f_1)-\mathcal{P}(f_2)||\nonumber\\
  &=&\sup_{(x,i)\in\mathbb{R}^+\times \mathcal{S}}|V_{f_1}(x,i)-V_{f_2}(x,i)|
  =\sup_{(x,i)\in\mathbb{R}^+\times \mathcal{S}}\left|\sup_{\pi\in\Pi}  R_{f_1,\pi}(x,i)-\sup_{\pi\in\Pi}R_{f_2,\pi}(x,i)\right|\nonumber\\
   &\le&\sup_{(x,i)\in\mathbb{R}^+\times \mathcal{S}}\sup_{\pi\in\Pi}  |R_{f_1,\pi}(x,i)-R_{f_2,\pi}(x,i)|
\sup_{(x,i)\in\mathbb{R}^+\times E}E_{(x,i)}\left[e^{-\Lambda_{\sigma_1}}||
    f_1-f_2||\right]\nonumber\\
  &=&||
  f_1-f_2||\int_0^\infty q_ie^{-q_i t}e^{-\delta_i t}\dif t=\max_{i\in E}\frac{q_i}{q_i+\delta_i}||f_1-f_2||,\label{jth2}
\end{eqnarray}
where the last inequality follows by \eqref{4214-2} and the last equality follows by noting that $\sigma_1$ is exponentially distributed with mean $\frac1{q_i}$. Therefore,   $\mathcal{P}$ is a contraction on the space
$(\mathcal{D},||\cdot||)$.

Note that for any $f\in\mathcal{C}$ and $i\in\mathcal{S}$, $f(\cdot,i)$ is non-decreasing. Hence, it follows by \eqref{4214-2} and \eqref{6214-1} that the operator $\mathcal{P}$ is non-decreasing. Consider two functions $g_1,g_2$ defined by $g_1(x,i)=0$ and $g_2(x,i)=\frac{\overline{l}}{\underline{\delta}(1+d)}$. It is not hard to verify that $g_1,g_2\in\mathcal{D}$ and $g_1\le V\le g_2$. Hence, $\mathcal{P}(g_1)\le\mathcal{P}(V)\le \mathcal{P}(g_2)$. Note that by \eqref{60214-3} $\mathcal{P}(V)=V$. Hence, $\mathcal{P}(g_1)\le V\le \mathcal{P}(g_2)$. Apply the operator $\mathcal{P}$ once again, we have $\mathcal{P}^2(g_1)\le V\le \mathcal{P}^2(g_2)$.     By repeating this $n-2$ more times, we obtain $\mathcal{P}^n(g_1)\le V\le \mathcal{P}^n(g_2)$.  As a result, $\lim_{n\rightarrow}\mathcal{P}^n(g_1)\le V\le \lim_{n\rightarrow \infty}\mathcal{P}^n(g_2)$. Since $\mathcal{P}$ is a contraction on the complete space $(\mathcal{D},||\cdot||)$, there is a unique fixed point in $\mathcal{D}$ and  is identical to both $\lim_{n\rightarrow \infty}\mathcal{P}^n(g_1)$ and $\lim_{n\rightarrow \infty}\mathcal{P}^n(g_2)$. Consequently, $\lim_{n\rightarrow \infty}\mathcal{P}^n(g_2)=V=\lim_{n\rightarrow \infty}\mathcal{P}^n(g_2)$.  As a result, $V\in\mathcal{D}$.
\noindent (ii) The results follow immediately by (i) and Theorem \ref{23513-11}.
\hfill $\square$

\noindent {\bf Proof of Theorem \ref{thmj1}}
Since, $ b^V_i< \infty$ for all $i\in \mathcal{S}$,
we can define an operator $\mathcal{Q}$ by
\begin{eqnarray}
  \mathcal{Q}(f)(x,i)=R_{f,\pi^{0,b^V_i}}(x,i)\ \mbox{
    for $f\in\mathcal{C}$, $x\ge 0$,  and $i\in \mathcal{S}$}.\label{jth6}
\end{eqnarray}
  The function $R_{f,\pi^{0,b^V_i}}$ is obviously nonnegative according to its definition. It follows by Lemma \ref{remff9} that $R_{f,\pi^{0,b^V_i}}\le V_f\le \frac{\bar{l}}{{\underline{\delta}(1+d)}}$ and by Corollary \ref{the2} that the function $R_{f,\pi^{0,b^V_i}}(\cdot,i)$ is increasing. Therefore, $R_{f,\pi^{0,b^V_i}}\in\mathcal{C}$. Then by \eqref{jth6} we know $\mathcal{Q}(f)\in\mathcal{C}$.
  It follows by \eqref{4214-2} that
\begin{eqnarray*}
  ||\mathcal{Q}(f_1)-\mathcal{Q}(f_2)||
   &=&\sup_{(x,i)\in\mathbb{R}^+\times \mathcal{S}}|R_{f_1,\pi^{0,b^V_i}}(x,i)-R_{f_2,\pi^{0,b^V_i}}(x,i)|\nonumber\\
     &\le&\sup_{(x,i)\in\mathbb{R}^+\times E}E_{(x,i)}\left[e^{-\Lambda_{\sigma_1}}||
    f_1-f_2||\right]\nonumber\\
  &=&||
  f_1-f_2||\int_0^\infty q_ie^{-q_i t}e^{-\delta_i t}\dif t=\max_{i\in E}\frac{q_i}{q_i+\delta_i}||f_1-f_2||.
\end{eqnarray*}
Consequently, $\mathcal{Q}$ is a contraction on $(\mathcal{C},||\cdot||)$. Hence, there is a unique fixed point of $\mathcal{Q}$ on $(\mathcal{C},||\cdot||)$.
Note by \eqref{jth6} we have
$ \mathcal{Q}(V)(x,i)=R_{V,\pi^{0,b^V_i}}(x,i)=V(x,i)$, where the last equality follows by Theorem \ref{6214-4} (ii). Therefore, $V$ is a fixed point.
By  \eqref{jth6}  and noticing that $\pi^{0,b_i^V}$ and $\pi^*$ are identical before $\sigma_1$,  we have
 \begin{align}\mathcal{Q}(R_{\pi^*})(x,i)&=R_{R_{\pi^*},\pi^{0,b_i^V}}(x,i)
  \\
  &=\E_{(x,i)}\bigg[\int_0^{  \sigma_1} e^{-\Lambda_t}\frac{l_t^*}{1+d}\dif t-\int_0^{  \sigma_1} e^{-\Lambda_t}\frac1{1-c}\dif C_t^*\nonumber\\
  &+e^{-\Lambda_{\sigma_1}}R_{\pi^*}(X^{\pi^*}_{\sigma_1}, \xi_{\sigma_1})
    \bigg],\ x\ge 0,i\in\mathcal{S},\label{4214-21}
\end{align}
where the last equality follows by \eqref{4214-2}.
It is not hard to see that the process $(X^{\pi^*},J)$ is a Markov process. Hence, it follows by the Markov property that
\begin{align}R_{\pi^*}(x,i)
  &=\E_{(x,i)}\bigg[\int_0^{  \sigma_1} e^{-\Lambda_t}\frac{l_t^*}{1+d}\dif t-\int_0^{  \sigma_1} e^{-\Lambda_t}\frac1{1-c}\dif C_t^*\nonumber\\
  &+e^{-\Lambda_{\sigma_1}}R_{\pi^*}(X^{\pi^*}_{\sigma_1}, \xi_{\sigma_1})
    \bigg],\ x\ge 0,i\in\mathcal{S}.\label{4214-20}
\end{align}
Combining \eqref{4214-21} and \eqref{4214-20} we obtain
$
\mathcal{Q}(R_{\pi^*})(x,i)=R_{\pi^*}(x,i), \ x\ge 0,i\in\mathcal{S}.$
Therefore, $R_{\pi^*}$ is also a fixed point. As there is a unique fixed point, we conclude $V=R_{\pi^*}$.\hfill $\square$\\

\noindent {\bf Acknowledgements}
  This work was supported by the University of New South Wales Australian Business School Special Research Grants.

%

\begin{thebibliography}{}

\bibitem[B{\"a}uerle, 2004]{Bauerle2004}
B{\"a}uerle, N. (2004).
\newblock Approximation of optimal reinsurance and dividend payout policies.
\newblock {\em Mathematical Finance}, 14(1):99--113.

\bibitem[Cadenillas et~al., 2007]{CadenillasSarkarZapatero2007}
Cadenillas, A., Sarkar, S., and Zapatero, F. (2007).
\newblock Optimal dividend policy with mean-reverting cash reservoir.
\newblock {\em Mathematical Finance}, 17:81--109.

\bibitem[Dickson and Waters, 2004]{DicksonWaters2004}
Dickson, D.~C. and Waters, H.~R. (2004).
\newblock Some optimal dividends problems.
\newblock {\em ASTIN BULLETIN}, 34(1):49--74.

\bibitem[Fleming and Soner, 1993]{FlemingSoner1993}
Fleming, W.~H. and Soner, H.~M. (1993).
\newblock {\em Controlled {M}arkov processes and viscosity solutions}.
\newblock Applications of Mathematics. Springer-Verlag, New York.

\bibitem[He and Liang, 2008]{HeLiang2008}
He, L. and Liang, Z. (2008).
\newblock Optimal financing and dividend control of the insurance company with
  proportional reinsurance policy.
\newblock {\em Insurance: Mathematics and Economics}, 42(3):976 -- 983.

\bibitem[H{\o}jgaard and Taksar, 2001]{HojgaardTaksar2001}
H{\o}jgaard, B. and Taksar, M. (2001).
\newblock Optimal risk control for a large corporation in the presence of
  returns on investments.
\newblock {\em Finance and Stochastics}, 5(4):527--547.

\bibitem[Ikeda and Watanabe, 1977]{IkedaWatanabe1977}
Ikeda, N. and Watanabe, S. (1977).
\newblock A comparison theorem for solutions of stochastic differential
  equations and its applications.
\newblock {\em Osaka J. Math.}, 14(3):619--633.

\bibitem[Jiang and Pistorius, 2012]{JiangPistorius2012}
Jiang, Z. and Pistorius, M. (2012).
\newblock Optimal dividend distribution under {M}arkov regime switching.
\newblock {\em Finance and Stochastics}, 16(3):449--476.

\bibitem[Krylov, 1996]{Krylov1996}
Krylov, N.~V. (1996).
\newblock {\em Lectures on Elliptic and Parabolic Equations in H\"older
  Spaces}.
\newblock The American Mathematical Society.

\bibitem[L{\o}kka and Zervos, 2008]{LokkaZervos2008}
L{\o}kka, A. and Zervos, M. (2008).
\newblock Optimal dividend and issuance of equity policies in the presence of
  proportional costs.
\newblock {\em Insurance: Mathematics and Economics}, 42:954--961.

\bibitem[Paulsen, 2008]{Paulsen2008}
Paulsen, J. (2008).
\newblock Optimal dividend payments and reinvestments of diffusion processes
  with both fixed and proportional costs.
\newblock {\em SIAM Journal on Control and Optimization}, 47(5):2201--2226.

\bibitem[Scheer and Schmidli, 2011]{ScheerSchmidli2011}
Scheer, N. and Schmidli, H. (2011).
\newblock Optimal dividend strategies in a {C}ram\'er-{L}undberg model with
  capital injections and administration costs.

\bibitem[Shreve et~al., 1984]{ShreveLehoczkyGaver1984}
Shreve, S.~E., Lehoczky, J.~P., and Gaver, D.~P. (1984).
\newblock Optimal consumption for general diffusions with absorbing and
  reflecting barriers.
\newblock {\em SIAM Journal on Control and Optimization}, 22(1):55--75.

\bibitem[Sotomayor and Cadenillas, 2011]{SotomayorCadenillas2011}
Sotomayor, L.~R. and Cadenillas, A. (2011).
\newblock Classical and singular stochastic control for the optimal dividend
  policy when there is regime switching.
\newblock {\em Insurance: Mathematics and Economics}, 48(3):344 -- 354.

\bibitem[Taksar, 2000]{Taksar2000}
Taksar, M.~I. (2000).
\newblock Optimal risk and dividend distribution control models for insurance
  company.
\newblock {\em Mathematical Methods of Operations Research}, 51:1--42.

\bibitem[Yao et~al., 2011]{YaoYangWang2011}
Yao, D., Yang, H., and Wang, R. (2011).
\newblock Optimal dividend and capital injection problem in the dual model with
  proportional and fixed transaction costs.
\newblock {\em European Journal of Operational Research}, 211(3):568 -- 576.

\bibitem[Zhu, 2014a]{Zhu2014c}
Zhu, J. (2014a).
\newblock Dividend optimization for a regime-switching diffusion model with
  restricted dividend rates.
\newblock {\em ASTIN Bulletin}, 44:459--494.

\bibitem[Zhu, 2014b]{Zhu2014b}
Zhu, J. (2014b).
\newblock Dividend optimization for general diffusions with restricted dividend
  payment rates.
\newblock {\em Scandinavian Actuarial Journal, in Press}.

\bibitem[Zhu and Chen, 2013]{ZhuChen2013}
Zhu, J. and Chen, F. (2013).
\newblock Dividend optimization for regime-switching general diffusions.
\newblock {\em Insurance: Mathematics and Economics}, 53(2):439 -- 456.

\end{thebibliography}

\end{document}